\documentclass[reqno]{amsart}
\usepackage{amsrefs}
\usepackage[shortlabels]{enumitem}
\usepackage{hyperref}
\usepackage{cleveref}
\usepackage{tikz-cd}
\usepackage{bm}
\usepackage{physics}
\usepackage{geometry}

\theoremstyle{plain}
\newtheorem{thm}{Theorem}

\newtheorem{prop}[thm]{Proposition}
\newtheorem{lem}[thm]{Lemma}

\theoremstyle{definition}
\newtheorem{defn}[thm]{Definition}

\theoremstyle{remark}
\newtheorem{rk}[thm]{Remark}

\AddToHook{env/prop/begin}{\crefalias{thm}{prop}}
\AddToHook{env/lem/begin}{\crefalias{thm}{lem}}
\AddToHook{env/rk/begin}{\crefalias{thm}{rk}}

\crefname{prty}{property}{properties}
\crefname{defn}{Definition}{Definitions}
\crefname{thm}{Theorem}{Theorems}
\crefname{cor}{Corollary}{Corollaries}
\crefname{lem}{Lemma}{Lemmas}
\crefname{rk}{remark}{remarks}
\crefname{prop}{Proposition}{Propositions}

\DeclareMathOperator{\Id}{Id}

\DeclareMathOperator{\spn}{span}
\DeclareMathOperator{\ran}{ran}
\DeclareMathOperator{\supp}{supp}

\DeclareMathOperator{\proj}{proj}
\newcommand{\vertiii}[1]{{\left\vert\kern-0.25ex\left\vert\kern-0.25ex\left\vert #1 \right\vert\kern-0.25ex\right\vert\kern-0.25ex\right\vert}}
\newcommand{\fk}{{F_w(\mathcal H_\tau)}}

\newcommand{\hfk}{\hat F_w(\mathcal H_\tau)}

\title[Koopman, transfer operator techniques and quantum theory]{Koopman and transfer operator techniques from the perspective of quantum theory}
\author[Giannakis]{Dimitrios Giannakis}
\address{Department of Mathematics, Dartmouth College, Hanover, NH 03755, USA.}
\email{dimitrios.giannakis@dartmouth.edu}
\author[Montgomery]{Michael Montgomery}
\address{Department of Mathematics, Dartmouth College, Hanover, NH 03755, USA.}
\email{michael.r.montgomery@dartmouth.edu}

\begin{document}

\begin{abstract}
	The study of mathematical connections between operator-theoretic formulations of classical dynamics and quantum mechanics began at least as early as the 1930s in work of Koopman and von Neumann and was developed in later decades by many authors, often independently, into a framework now broadly known as Koopman--von Neumann representation of classical dynamics.
	This article surveys aspects of this framework for measure-preserving ergodic dynamical systems and connects it with recent approximation techniques for Koopman and transfer operators that are amenable to data-driven numerical implementation.
	In broad terms, these methods are based on representations of (i) classical observables as elements of an algebra of operators acting on a Hilbert space; and (ii) classical probability measures as elements of the state space of that algebra, with lifted versions of the Koopman and transfer operators inducing dynamical evolution of observables and states, respectively.
	A common theme underlying the techniques surveyed here is the use of reproducing kernel Hilbert spaces with coalgebra structure (so-called ``reproducing kernel Hilbert algebras'') that aids the quantum representation of classical objects, as well as the use of Fock spaces to build approximation schemes with high expressivity and structure preservation properties (notably, preservation of positivity and multiplicativity of composition operators). 
	Applications to quantum algorithms for approximating the Koopman evolution of observables in systems with pure point spectra are also discussed. 
\end{abstract}

\maketitle

\section{Introduction}
\label{sec:intro}

Let $\Phi\colon X \to X$ be an invertible, measure-preserving transformation of a probability space $(X, \Sigma, \mu)$.
A classical operator-theoretic approach, dating back to work of Koopman and von Neumann \cites{Koopman31,KoopmanVonNeumann32}, is to study this transformation through its induced (linear) action on the spaces of observables $L^p(\mu)$.
On these spaces the dynamics acts by means of the Koopman (or composition) operator,
\begin{displaymath}
	U \colon L^p(\mu) \to L^p(\mu), \quad U f = f \circ \Phi,
\end{displaymath}
which is an isometric isomorphism of $L^p(\mu)$, $p \in [1, \infty]$.
For $p \in (1,\infty]$, the Koopman operator is the dual of the Perron--Frobenius (or transfer) operator,
\begin{displaymath}
	P \colon L^q(\mu) \to L^q(\mu), \quad P f = f \circ \Phi^{-1},
\end{displaymath}
where $\frac{1}{p} + \frac{1}{q} = 1$.
Letting $\langle \cdot, \cdot \rangle_\mu$ denote the natural pairing between $L^p(\mu)$ and $L^q(\mu)$, we have
\begin{displaymath}
	\langle P w, f\rangle_\mu = \langle  w, U f\rangle_\mu, \quad \forall f \in L^p(\mu), \quad \forall w \in L^q(\mu),
\end{displaymath}
so there is a duality between the action of $\Phi$ on observables in $L^p(\mu)$ and absolutely continuous measures $\nu \ll \mu$ with densities $w = \frac{d\nu}{d\mu} \in L^q(\mu)$.

The operator-theoretic perspective lies at the core of many foundational results in ergodic theory.
For instance, $\Phi\colon X \to X$ is ergodic with respect to $\mu$ iff the eigenvalue of $U$ at 1 is simple, and it is weak-mixing iff $U$ has no eigenvalues other than 1 \cite{Walters81}.
Juxtaposed with weak-mixing systems are systems with pure point spectra; that is, measure-preserving transformations for which the union of eigenspaces of the Koopman operator $U$ is dense in $L^p(\mu)$, $p \in [1, \infty)$.
The Halmos--von Neumann theorem \cite{HalmosVonNeumann42} asserts that the point spectrum of $U$, denoted as $\sigma_p(U)$, completely characterizes such systems.
In more detail, two measure-preserving, ergodic, pure-point-spectrum systems on standard probability spaces, with Koopman operators $U$ and $\tilde U$, are measure-theoretically isomorphic iff $\sigma_p(U) = \sigma_p(\tilde U)$.
Moreover, every such system is isomorphic to a rotation system on a compact abelian group, and the point spectrum $\sigma_p(U)$ is a subgroup of $\mathbb T^1$.
Under continuous-time dynamics, a measure-preserving flow $\Phi^t \colon X \to X$, $t \in \mathbb R$, of appropriate regularity induces a strongly continuous unitary Koopman group $\{U^t \}_{t \in \mathbb R}$ on $L^2(\mu)$, $U^t f = f \circ \Phi^t$.
By Stone's theorem on one-parameter unitary groups \cite{Stone32}, this group is completely characterized by its generator, $V(D(V)) \to L^2(\mu)$; a skew-adjoint operator defined on a dense subspace $D(V) \subseteq L^2(\mu)$ via the norm limit $V f = \lim_{t\to 0} (U^t f - f) /t$ and giving $U^t = e^{tV}$ through the Borel functional calculus.
A fundamental property of the generator $V$ is that it acts as a derivation on the algebra $D(V) \cap L^\infty(\mu)$, satisfying the Leibniz rule,
\begin{equation}
	\label{eq:leibniz}
	V(fg) = (V f) g + f (V g), \quad \forall f, g \in D(V) \cap L^\infty(\mu).
\end{equation}
In fact, \eqref{eq:leibniz} characterizes which skew-adjoint operators on $L^2(\mu)$ are generators of unitary Koopman groups of composition operators \cite{TerElstLemanczyk17}.
See \cites{Halmos56,Baladi00,EisnerEtAl15} for expositions of these and other operator-theoretic results in ergodic theory.

In addition to their profound role in analytical studies, in  recent years, operator methods have seen widespread use in numerical techniques for data analysis, prediction, and control \cites{BruntonEtAl22,OttoRowley21,Colbrook24,MauroyEtAl20}.
In broad terms, these methods build approximations of the Koopman and/or transfer operator given samples of one or more observables taken along a dynamical trajectory or a collection of trajectories, with the type of approximation depending on the task at hand.
For instance, in prediction problems it is natural to consider approximations $\hat U$ of $U$ in the strong operator topology (having in mind that $\hat U f \approx U f$ is a prediction of a given observable $f$ at one timestep ``in the future''), leading to a fruitful interplay with the theory of statistical learning \cite{CuckerSmale01}.

In feature extraction applications, the focus is oftentimes on approximating spectral data such as elements of the spectra of $U$, $P$, and/or the generator $V$ and the corresponding spectral measures.
These methods typically require a careful choice of function space and stronger forms of approximation than prediction problems.
For example, early work from the 1990s \cite{Froyland97} used the Ulam method \cite{Ulam64} to approximate the transfer operator $P$ of uniformly expanding maps in spaces of H\"older-continuous functions, where $P$ is quasicompact, the goal being to identify observables that simultaneously have high regularity and slow correlation decay.
In \cite{DellnitzJunge99} the eigenvalues and eigenvectors of the transfer operator, approximated again with the Ulam method, were employed to identify near-invariant and near-periodic sets under stochastically perturbed dynamics with continuous transition functions.
On the Koopman operator side, data-driven techniques were initiated by \cite{MezicBanaszuk99,MezicBanaszuk04,Mezic05} who used harmonic averaging to approximate the point spectrum in $L^2$ and associated invariant and periodic phase space partitions.
The paper \cite{RowleyEtAl09} interpreted the dynamic mode decomposition (DMD) technique \cite{SchmidSesterhenn08,Schmid10} as a method for approximating Koopman eigenvalues and associated spatial projections now known as Koopman modes, opening up the way to Koopman spectral computations from high-dimensional snapshot data.
Algorithmically, DMD is closely related to principal oscillation pattern analysis and linear inverse modeling, which have been popular mode decomposition methods in the geophysical sciences since the 1980s; e.g., \cite{Hasselmann88,Penland89}.

More recently, a problem that has received significant attention is how to consistently represent the continuous spectra of Koopman or transfer operators on $L^2$ associated with mixing dynamics; e.g., \cites{KordaEtAl20,DasEtAl21,Colbrook23,ColbrookTownsend24,GiannakisValva25}.
Other approaches have emphasized the use of function spaces of higher regularity such as spaces of analytic functions  \cite{Wormell25,SlipantschukEtAl20,Mezic20}, fractional Sobolev spaces \cite{FroylandEtAl14b}, reproducing kernel Hilbert spaces (RKHSs) \cite{Kawahara16,DasGiannakis20,RosenfeldEtAl22}, reproducing kernel Banach spaces \cite{IkedaEtAl22}, and anisotropic Banach spaces adapted to the dynamics \cite{BlankEtAl02}, in order to detect isolated eigenvalues and resonances.
For further discussion on the literature on spectral approximation techniques for Koopman and transfer operators and their applications we refer the reader to the review papers \cite{OttoRowley21,Colbrook24,BruntonEtAl22}, as well as other chapters in this collection.

\subsection{Operator-theoretic techniques and quantum mechanics}
\label{sec:qm_overview}

Besides providing powerful tools for analyzing nonlinear state space dynamics, the operator-theoretic formalism also provides a connection between classical and quantum theories.

From an algebraic standpoint, one associates to a classical dynamical system a unital, abelian algebra of observables on which the dynamics acts, through the Koopman operator, as an algebra morphism.
In the setting of measure-preserving, invertible transformations of the probability space $(X, \Sigma, \mu)$ considered above a natural choice is $\mathfrak A := L^\infty(\mu)$ which is a von Neumann algebra \cite{Takesaki01} (i.e., a commutative $C^*$-algebra with a Banach space predual, $\mathfrak A_* := L^1(\mu)$) with respect to pointwise function multiplication and complex conjugation.
The Koopman operator acts on this algebra as a $^*$-isomorphism,
\begin{displaymath}
	U(fg) = (U f)(U g), \quad U(f^*) = (U f)^*, \quad \forall f, g \in \mathfrak A,
\end{displaymath}
with the Leibniz rule~\eqref{eq:leibniz} being an analog of this property in continuous time.

In quantum theory, on the other hand, the algebra of observables, $\mathfrak B$, is non-abelian---typically, it is chosen as a unital $C^*$-algebra $\mathfrak B \subseteq B(H)$ for a Hilbert space $H$ and space of bounded operators $B(H)$ equipped with operator composition and adjoint \cite{Pillet06}.
Considering again a measure-preserving, invertible transformation $\Phi$ of $(X, \Sigma, \mu)$, it is natural to choose the Hilbert space $H=L^2(\mu)$ (on which the Koopman and transfer operators act as unitaries, $U^* \equiv U^{-1} \equiv P$), and set $\mathfrak B = B(H)$.
With this choice, $\mathfrak B$ is a von Neumann algebra with predual $\mathfrak B_* = B_1(H)$ (the space of trace-class operators on $H$ equipped with the trace norm).
The unitary Koopman operator on $H$ has an adjoint action $\mathcal U \colon \mathfrak B \to \mathfrak B$, $\mathcal U A = U A U^*$, which is again a $^*$-isomorphism,
\begin{displaymath}
	\mathcal U(A B) = (\mathcal U A) (\mathcal U B), \quad \mathcal U(A^*) = (\mathcal U A)^*, \quad \forall A, B \in \mathfrak B.
\end{displaymath}

Recall now that $\mathfrak A = L^\infty(\mu)$ admits a representation $\pi \colon \mathfrak A \to \mathfrak B $ into $\mathfrak B = B(L^2(\mu))$ that assigns to $f \in \mathfrak A$ the multiplication operator $\pi f \in \mathfrak B$, where $(\pi f) g = f g$.
One readily verifies that $\pi$ is a faithful $^*$-representation (that is, $\pi f = \pi \tilde f$ iff $f = \tilde f$ and $(\pi f)^* = \pi(f^*)$) that is compatible with Koopman evolution, $\pi \circ U = \mathcal U \circ \pi$.
The latter relationship can be depicted as a commuting diagram,
\begin{equation}
	\label{eq:commut_koopman}
	\begin{tikzcd}
		\mathfrak A \ar[r,"U"] \ar[d,"\pi",swap] & \mathfrak A \ar[d,"\pi"] \\
		\mathfrak B \ar[r,"\mathcal U"] & \mathfrak B
	\end{tikzcd},
\end{equation}
which illustrates that one can think of $\pi$ as an embedding of dynamics of classical observables in $\mathfrak A$ into dynamics of quantum observables in $\mathfrak B$.

A dual picture to the above follows by considering the evolution of states (i.e., positive, unital functionals) of the algebras $\mathfrak A$ and $\mathfrak B$.
With $\mathfrak A = L^\infty(\mu)$, let $ S_*(\mathfrak A) $ be the set of probability densities in $\mathfrak A_* = L^1(\mu)$; that is, $S_*(\mathfrak A)$ consists of positive elements $\sigma \in L^1(\mu)$ satisfying $\int_X \sigma \, d\mu = 1$.
Every $\sigma \in S_*(\mathfrak A)$ induces a state $\mathbb E_\sigma \colon \mathfrak A \to \mathbb C$ that implements statistical expectation, $\mathbb E_\sigma f := \int_X \sigma f \, d\mu$.
Such states are called normal states.
One can directly verify that $\mathfrak A_*$ is invariant under the action of the transfer operator $P$ on $L^1(\mu)$, and the duality relationship $\mathbb E_{P \sigma} f = \mathbb E_\sigma(U f)$ holds for every $f \in \mathfrak A$ and $\sigma \in \mathfrak A_*$.

In the quantum setting of $\mathfrak B = B(H)$, $H=L^2(\mu)$, normal states are induced by density operators; that is, positive, trace-class operators in $\mathfrak B_* = B_1(H)$ with unit trace.
Denote the set of all such operators as $S_*(\mathfrak B)$.
Every $\rho \in S_*(\mathfrak B)$ induces a normal state $\mathbb E_\rho \colon \mathfrak B \to \mathbb C$ defined as $\mathbb E_\rho A = \tr(\rho A)$.
Analogously to $\mathbb E_\sigma f$ from the classical setting, the quantity $\mathbb E_\rho A$ is the expectation of quantum observable $A$ with respect to the state induced by $\rho$.
Moreover, defining $\mathcal P \colon B_1(H) \to B_1(H)$ as $\mathcal P A = U^* A U$, we have the duality relationship $\mathbb E_{\mathcal P\rho} A = \mathbb E_\rho(\mathcal U A)$ for every quantum observable $A \in \mathfrak B$ and density operator $\rho \in S_*(\mathfrak B)$, again analogously to the classical case.
In particular, one can think of $\mathcal P$ as a quantum version of the transfer operator $P$.
In quantum mechanics, the evolution of observables by $\mathcal U$ and the evolution of states by $\mathcal P$ are oftentimes referred to as the Heisenberg and Schr\"odinger pictures, respectively.

To complete our classical--quantum correspondence, we employ an embedding of classical probability densities in $S_*(\mathfrak A)$ into density operators in $S_*(\mathfrak B)$ that plays a dual role to the embedding $\pi$ of observables, with the key difference that the embedding of densities will be nonlinear.
Specifically, we define $\Gamma\colon S_*(\mathfrak A) \to S_*(\mathfrak B)$ as $\Gamma(\sigma) = \langle \sigma^{1/2}, \cdot\rangle_H \sigma^{1/2}$, and it follows from the fact that $\sigma$ is a probability density that $\Gamma(\sigma)$ is a rank-1 density operator that projects along the unit vector $\sigma^{1/2} \in H$.
The corresponding state $\mathbb E_{\Gamma (\sigma)}$ is known as a vector state of $\mathfrak B$ and it is an example of a pure state (i.e., an extremal point in the set of states of $\mathfrak B$, equipped with the weak-$^*$ topology).

The map $\Gamma$ is an example of embeddings of probability densities into quantum states, now broadly known as Koopman--von Neumann embeddings \cite{Mauro02,Barandes25}.
Since early references \cite{Schonberg52,Schonberg53,DellaRiciaWiener66}, first appearing in the 1950s, the Koopman--von Neumann approach has found applications in diverse mathematical and applied areas, including semigroup theory, statistical mechanics, quantum simulation, plasma dynamics, kinetic theory, and hybrid classical--quantum systems; e.g., \cite{Joseph20,LinEtAl22,JosephEtAl23,StenglEtAl24,NovikauJoseph25}.
The formalism has also been extended and applied in areas including semiclassical \cite{BondarEtAl12,BondarEtAl19,JosephEtAl23} and relativistic \cite{Giannakis21b,Mezic23} dynamics.

Returning to the measure-preserving examples at hand, the map $\Gamma$ intertwines the classical and quantum transfer operators, $\Gamma \circ P = \mathcal P \circ \Gamma $, and satisfies an analogous commuting diagram to~\eqref{eq:commut_koopman},
\begin{equation}
	\label{eq:commut_transfer}
	\begin{tikzcd}
		S_*(\mathfrak A) \ar[r,"P"] \ar[d,"\Gamma",swap] & S_*(\mathfrak A) \ar[d,"\Gamma"] \\
		S_*(\mathfrak B) \ar[r,"\mathcal P"] & S_*(\mathfrak B)
	\end{tikzcd}.
\end{equation}
Thus, similarly to the setting of observables, we have a dynamically consistent embedding of classical probability densities (representing normal classical statistical states) into quantum density operators (representing normal quantum states).
In summary, the correspondences between classical--quantum evolution and Heisenberg--Schr\"odinger picture described in this subsection can be expressed as
\begin{equation}
	\label{eq:consistency}
	\mathbb E_\sigma(U f) = \mathbb E_{P \sigma} f = \mathbb E_{\Gamma(\sigma)}(\mathcal U (\pi f)) = \mathbb E_{\mathcal P \Gamma(\sigma)}(\pi f), \quad \forall f \in \mathfrak A, \quad \forall \sigma \in S_*(\mathfrak A).
\end{equation}

\subsection{Objectives and plan of the paper}

This expository article surveys Koopman and transfer operator techniques that utilize, in some level, embeddings of classical dynamics into non-commutative operator algebras, similar to \eqref{eq:commut_koopman} and \eqref{eq:commut_transfer}.
The material presented is mainly taken from \cite{Giannakis19b,FreemanEtAl23,DasGiannakis23,DasEtAl23,GiannakisEtAl24,GiannakisEtAl25,GiannakisEtAl22}.
Using these methods as examples, our goal is to highlight properties of non-commutative spaces of operators that are useful in analysis and approximation techniques for classical dynamical systems.
In particular, we focus on (i) positivity preservation of observables and algebra states under finite-dimensional projection; and (ii) recovery of the Leibniz rule for the generator after spectral regularization.
We should note that a general review of methods for quantum simulation of classical dynamics lies outside the scope our discussion.
Methods based on truncated Carleman linearization \cite{LiuEtAl21,WuEtAl25,CostaEtAl25}, matrix product states \cite{Oseledets10,Khoromskij11,YeLoureiro24}, Lie--Trotter--Suzuki approximation \cite{LloydEtAl20}, Schr\"odingerization \cite{JinEtAl24}, linear combination of unitaries \cite{AnEtAl23,BharadwajSreenivasan25}, quantum phase estimation \cite{LokareEtAl24}, quantum reservoir computers \cite{PfefferEtAl22}, and iterated quantum measurements \cite{AndressEtAl24} are recent examples from the literature on this area that the interested reader may wish to refer to.
Another important topic not covered in this article is that of numerical approximation.
While we do not discuss these aspects here, all of the methods presented are amenable to data-driven approximation schemes with convergence guarantees in appropriate large-data limits.
Details on these techniques, which are largely based on kernel methods for machine learning and RKHS theory, can be found in the references cited above.

The remainder of the paper is organized as follows.
We begin in \cref{sec:qmda} by surveying a quantum mechanical data assimilation (QMDA) scheme \cite{Giannakis19b,FreemanEtAl23} that represents uncertain knowledge about the state of a partially observed dynamical system by means of quantum density operators (as opposed to probability densities in the classical case).
In \cref{sec:rkha}, we review properties of a class of function spaces, called reproducing kernel Hilbert algebras (RKHAs) \cite{DasGiannakis23,DasEtAl23,GiannakisMontgomery25}, which are RKHSs with additional algebra and coalgebra structure that is useful for building quantum mechanical embeddings.
In \cref{sec:overview_fock_space} we describe methods that leverage the properties of RKHAs to lift the Koopman/transfer operators of measure-preserving systems to evolution operators acting on a Fock space associated with a many-body quantum system \cite{GiannakisEtAl24,GiannakisEtAl25}.
In \cref{sec:qc}, we review an RKHA-based technique \cite{GiannakisEtAl22} for simulating continuous-time systems with pure point spectra on quantum computers.

Unless otherwise stated, in what follows the state space $X$ will be assumed to be a compact metrizable Hausdorff space, the dynamics $\Phi\colon X \to X$ (or $\Phi^t \colon X \to X$ in continuous time) will be assumed to be continuous and invertible, and $\mu$ will be assumed to be a dynamically invariant Borel probability measure.
A number of these assumptions can be relaxed but we make them here to avoid technicalities.
We will let $\iota\colon C(X) \to L^2(\mu)$ denote the (continuous) restriction map from the Banach space of continuous functions on $X$ equipped with the max norm, to $L^2(\mu)$.

\section{Quantum mechanical data assimilation}
\label{sec:qmda}

Data assimilation is a framework for state estimation and prediction for partially observed dynamical systems \cites{MajdaHarlim12,LawEtAl15}.
Its sequential formulation, also known as filtering, is based on a predictor-corrector procedure, whereby a forward model is employed to evolve the probability distribution for the system state until a new observation is acquired, at which time that probability distribution is updated in an analysis step to a posterior distribution correcting for model error and/or uncertainty in the prior distribution.
Since the seminal work of Kalman \cite{Kalman60} on filtering (which utilizes Bayes' theorem for the analysis step under the assumption that all distributions are Gaussian), data assimilation has evolved to an indispensable tool in object tracking \cite{ThrunEtAl05}, weather forecasting \cite{Kalnay03}, and many other important applications involving complex systems \cite{LahozEtAl10}.
In this section, we survey a quantum mechanical formulation of data assimilation for measure-preserving dynamics \cite{Giannakis19b,FreemanEtAl23} that uses probability densities in $L^1(\mu)$ and density operators in $B_1(H)$, $H = L^2(\mu)$, to represent uncertain knowledge about the system state in the classical and quantum formulations, respectively.

\subsection{Quantum mechanical formulation of Bayesian filtering}
\label{sec:qm_filtering}

Suppose that the filtering probability density at timestep $n \in \mathbb N$ is $\sigma_n \in L^1(\mu)$.
Then, the forecast step in classical data assimilation amounts to pushing forward $\sigma_n$ under the transfer operator, i.e., $\tilde\sigma_{n+1} = P \sigma_n$.
This gives the prior density $\tilde\sigma_{n+1}$ at timestep $n+1$, before a new observation of the system is made.
We have already seen how to represent this step quantum mechanically in \cref{sec:qm_overview}, $\tilde\rho_{n+1} = \mathcal P \rho_n$, where $\rho_n \in B_1(H)$ is the quantum density operator at timestep $n$.
In particular, if we have
\begin{equation}
	\label{eq:qm_rep}
	\rho_n = \Gamma(\sigma_n),
\end{equation}
the quantum state $\mathbb E_{\rho_n}$ yields consistent statistical predictions about any classical observable $f \in L^\infty(\mu)$ in the sense of~\eqref{eq:consistency}.

Next, we consider a quantum mechanical formulation of the Bayesian analysis step.
Classically, we model an observation of the system at timestep $n+1$ as an event, i.e., a measurable set $\Omega_{n+1} \subseteq X$.
For example, for a system observed through an observable $h\colon X \to \mathbb R$ using a detector of resolution $\delta>0$, we have $\Omega_{n+1} = h^{-1}(I_{n+1})$ where $I_{n+1} \subset \mathbb R$ is an interval of length $\delta$ that contains the measurement of $h$ recorded by the detector.
Measurable observation maps $h\colon X \to Y$ taking values in a potentially higher-dimensional observation space $Y \subseteq \mathbb R^d$ can be treated similarly.

Let $\chi_{\Omega_{n+1}}\colon X \to \mathbb R$ be the characteristic function of $\Omega_{n+1}$.
Assuming that $\int_X\chi_{\Omega_{n+1}}\tilde\sigma_{n+1}\,d\mu > 0$, the Bayesian posterior $\sigma_{n+1} \in L^1(\mu)$ of $\tilde\sigma_{n+1}$ given $\Omega_{n+1}$ is
\begin{equation}
	\label{eq:cl_bayes}
	\sigma_{n+1} \equiv \tilde\sigma_{n+1}\rvert_{\Omega_{n+1}} = \frac{\tilde\sigma_{n+1}\chi_{\Omega_{n+1}}}{\int_X\tilde\sigma_{n+1}\chi_{\Omega_{n+1}}\,d\mu}.
\end{equation}
From an algebraic standpoint, the characteristic function $\chi_{\Omega_{n+1}}$ is a positive element of $L^\infty(\mu)$ that is bounded by the unit of the algebra, $0 \leq \chi_{\Omega_{n+1}}(x) \leq 1_X(x)$ for $\mu$-a.e.\ $x \in X$.

In quantum theory, the analogs of characteristic functions representing classical events are positive operators $e \in B(H)$ that are bounded above by the identity, $0 \leq e \leq \Id$.
Such operators are known as \emph{quantum effects}.
They form associated effect algebras that play a central role in formulations of quantum logic and quantum probability; e.g., \cite{Gudder07}.

For an effect $e \in B(H)$ and a density operator $\rho \in B_1(H)$ such that $\tr(\rho e) > 0$, the posterior density operator $\rho_n\lvert_e$ is given by
\begin{equation}
	\label{eq:qm_bayes}
	\rho\rvert_e = \frac{\sqrt e \rho \sqrt e}{\tr(\sqrt e \rho \sqrt e)},
\end{equation}
where $\sqrt{\cdot}$ denotes the square root of positive operators.
Intuitively, \eqref{eq:qm_bayes} can be thought of as a non-commutative Bayes' rule.

In our data assimilation framework, we can represent a classical event quantum mechanically as an effect that acts as a multiplication operator by the corresponding characteristic function---specifically, event $\Omega_{n+1}$ from~\eqref{eq:cl_bayes} is represented by effect $e_{n+1} \in B(H)$ where $e_{n+1} f = \chi_{\Omega_{n+1}} f$.
An application of the Bayes formula~\eqref{eq:qm_bayes} then yields the posterior density operator $\rho_{n+1} = \tilde\rho_{n+1}\rvert_{e_{n+1}}$ induced by event $\Omega_{n+1}$ observed at timestep $n+1$.
One readily verifies that under~\eqref{eq:qm_rep} $\rho_{n+1}$ is consistent with the classical posterior $\sigma_{n+1}$; that is, $\Gamma(\sigma_{n+1}) = \rho_{n+1}$.

A generalization of this approach is to employ a kernel function $\kappa\colon Y \times Y \to [0, 1]$ with associated feature map $\varphi\colon Y \to H$ given by $\varphi(y) = \kappa(y, h(\cdot))$.
Promoting $\varphi(y)$ to a multiplication operator then leads to an effect-valued feature map $\varepsilon\colon Y \to B(H)$, where $\varepsilon(y) f = \varphi(y) f$.
Letting $y_{n+1} = h(x_{n+1}) \in Y$ be the observation made at timestep $n+1$ with underlying classical state $x_{n+1} \in X$, the posterior density operator $\rho\rvert_{\varepsilon(y_{n+1})}$ is consistent with the posterior density $\sigma_{n+1}$
obtained by conditioning $\tilde\sigma_{n+1}$ by the likelihood function $\varepsilon(y_{n+1})$; i.e., $\sigma_{n+1} = \tilde\sigma_{n+1} \varphi(y_{n+1}) / Z$ with $Z = \int_X \varphi(y_{n+1}) \, d\mu$.

With either approach, iterating the quantum mechanical forecast and analysis steps described above yields a sequence $\rho_1, \rho_2, \ldots$ of density operators in $B_1(H)$ that is consistent (in the sense of~\eqref{eq:qm_rep} holding for every $n \in \mathbb N$) with the sequence $\sigma_1, \sigma_2, \ldots $ of probability densities in $L^1(\mu)$ obtained under classical Bayesian filtering for observations $y_n = h(x_n)$ induced by the classical dynamical trajectory $x_{n+1} = \Phi(x_n)$.

\subsection{Positivity-preserving projection}

Aside from special cases (e.g., the invariant measure $\mu$ being supported on a finite periodic orbit) the quantum mechanical system described in \cref{sec:qm_filtering} is infinite-dimensional.
To build a practical data assimilation algorithm we project this system to finite dimensions.

Let $\Pi_0, \Pi_1, \ldots$ be a sequence of finite-rank orthogonal projections on $H$ converging strongly to the identity, $\lim_{L\to\infty} \Pi_L f = f$ for every $f \in H$.
As a typical example, the family $\{\Pi_L\}_{L =0}^\infty$ can be constructed from an orthonormal basis $\{\phi_l\}_{l=0}^\infty$ of $H$; that is, $\Pi_L f = \sum_{l=0}^{L-1} \langle \phi_l, f \rangle_H \phi_l$.
We will henceforth assume such a choice of projections for concreteness.
Let $H_L \subset H$ be the range of $\Pi_L$.
The operator algebra $B(H_L)$ has dimension $L^2$ and is isomorphic with the algebra $\mathbb M_L(\mathbb C)$ of $L \times L$ complex matrices.
In particular, $A \in B(H_L)$ is represented by matrix $\bm A \in \mathbb M_L(\mathbb C)$ with elements $A_{ij} = \langle\phi_i, \phi_j \rangle_H$.
The algebra $B(H_L)$ can also be identified with the subalgebra of $B(H)$ consisting of all elements $A$ with $\ker A \supseteq H_L$ and $\ran A \subseteq H_L$.

Using the projections $\Pi_L$ we can compress operators $A \in B(H)$ to finite-rank operators in $B(H_L)$ through the linear map $\bm \Pi_L \colon B(H) \to B(H)$, where $A_L = \bm \Pi_L A := \Pi_L A \Pi_L$.
In particular, we have projected Koopman operators $U_L = \bm \Pi_L U$, projected multiplication operators (Toeplitz operators) $\pi_L f = \bm\Pi_L(\pi f)$, and effect-valued feature maps $\varepsilon_L = \bm \Pi_L \circ \varepsilon $.
We also have projected quantum density operators $\Gamma_L(\sigma) = \bm\Pi_L(\Gamma(\sigma)) / \tr(\bm\Pi_L(\Gamma(\sigma)))$ which are defined whenever $\bm\Pi_L(\Gamma(\sigma)) \neq 0$.
All of these maps converge as $L\to\infty$ in appropriate operator topologies; in particular, $U_L$, $\pi_L f$, and $\varepsilon_L(y)$ converge in the strong topology of $B(H)$, while $\Gamma_L(\sigma)$ converges in the norm topology of $B_1(H)$ for any probability density $\sigma \in L^1(\mu)$.

With these definitions, the finite-dimensional formulation of QMDA is essentially analogous to the infinite-dimensional scheme from \ref{sec:qmda}, with all operators in $B(H)$ replaced by their compressed versions in $B(H_L)$ obtained using $\bm \Pi_L$.
A key property of this projection is positivity preservation; that is, $\bm \Pi_L$ maps positive elements of $B(H)$ to positive elements of $B(H_L)$.
This is in contrast to orthogonal projections on $H$ which, in general, do not preserve positivity of functions (e.g., an $L^2$ orthogonal projection $\Pi_L f$ of a positive function $f$ may develop Gibbs oscillations to negative values).
In effect, by embedding classical data assimilation to the infinite-dimensional operator algebra $B(H)$ \emph{before} projecting to finite dimensions allows us to take advantage of properties of non-commutative spaces of operators (namely, positivity preservation under conjugation of operators) to build finite-dimensional approximations with structure preservation properties that are not available through approximations performed at the classical level.
In that regard, our approach aligns itself with the philosophy of avoiding discretization until the last possible step, which has been advocated in the context of Bayesian inverse problems \cite{Stuart10}.

Still, despite being positivity preserving and asymptotically consistent in the infinite dimension limit ($L\to\infty$) the finite-dimensional QMDA scheme has certain structural differences from its infinite-dimensional counterpart.
In some sense, the finite-dimensional scheme is more ``quantum mechanical" than the underlying infinite-dimensional system:
\begin{itemize}
	\item The multiplication operators $\{\pi f: f \in L^\infty(\mu) \} $ form a (maximal) abelian subalgebra of $B(H)$.
	      In contrast, the projected multiplication operators $\pi_L f$ generally do not commute, and generate a non-abelian algebra within $B(H_L)$.
	\item In general, the projected quantum sates $\Gamma_L(\sigma)$ are not the images of classical probability density functions under $\Gamma$.
	      That is, in general there is no probability density $\sigma_L \in L^1(\mu)$ such that $\Gamma_L(\sigma)$ is equal to $\langle\sigma_L^{1/2},\cdot\rangle_H\sigma_L^{1/2}$.
	      In fact, we have $\Gamma_L(\sigma) \propto \langle\psi_L, \cdot\rangle_H \psi_L$ with $\psi_L = \Pi_L \sigma^{1/2}$ playing the role of a quantum mechanical wavefunction that is not necessarily positive.
	\item The projected Koopman operator $U_L$ is not a composition operator.
	      Therefore, the mapping $f \mapsto U_L f$ is not interpretable as an evolution of observables under a classical flow, and nor is $\sigma \mapsto U_L^* \sigma$ interpretable as an evolution of densities.
	      In contrast, $A \mapsto U_L A U_L^*$ is well-defined quantum operation on $B(H)$ (i.e., a completely positive map with a non-trace-increasing predual; see \cite{Holevo01}) associated with an open quantum system.
	      If it happens that $H_L$ is an invariant subspace under $U$ then $U_L$ defines a quantum channel (i.e., a closed quantum system), but typically this will not be the case.
	      If preserving unitarity of the approximate Koopman dynamics on $H_L$ is of interest then the operator compression $U \mapsto \Pi_L U \Pi_L$ may be replaced by a unitary approximation obtained via one of the techniques available in the literature; e.g., \cite{DasEtAl21,Colbrook23,GiannakisValva24,GiannakisValva25}
\end{itemize}

\section{Reproducing kernel Hilbert algebras}
\label{sec:rkha}

We review the definition, properties, and constructions of RKHAs employed in subsequent sections of this chapter.
A detailed analysis of these objects can be found in \cite{GiannakisMontgomery25}.
For expositions of RKHS theory, see \cite{PaulsenRaghupathi16,SteinwartChristmann08}.

\subsection{Definitions and basic constructions}
\label{sec:rkha_defs}

Let $\mathcal H$ be an RKHS of complex-valued functions on a set $X$ with reproducing kernel $k\colon X \times X \to \mathbb C$ and inner product $\langle \cdot, \cdot \rangle_{\mathcal H}$.
We use $k_x := k(x, \cdot)$ to denote the kernel section at $x \in X$ and $\delta_x \colon \mathcal H \to \mathbb C$ to denote the corresponding pointwise evaluation functional, $\delta_x f = f(x) = \langle k_x, f\rangle_{\mathcal H}$.
We also let $\varphi \colon X \to \mathcal H$ denote the canonical feature map, $\varphi(x) = k_x$.

\begin{defn}[RKHA]
	An RKHS $\mathcal H$ on a set $X$ is a \emph{reproducing kernel Hilbert algebra (RKHA)} if $k_x \mapsto k_x \otimes k_x$, $x \in X$, extends to a bounded linear map $\Delta\colon \mathcal H \to \mathcal H \otimes \mathcal H$.
\end{defn}

Since
\begin{displaymath}
	\langle k_x, \Delta^*(f \otimes g)\rangle_{\mathcal H} = \langle \Delta k_x, f \otimes g \rangle_{\mathcal H \otimes \mathcal H} = \langle k_x \otimes k_x, f \otimes g\rangle_{\mathcal H \otimes \mathcal H} = \langle k_x, f\rangle_{\mathcal H} \langle k_x, g\rangle_{\mathcal H} = f(x)g(x),
\end{displaymath}
it follows that $\Delta^*$ is a bounded linear map that implements pointwise multiplication.
As a result, $\mathcal H$ is simultaneously a Hilbert function space and commutative algebra with respect to pointwise function multiplication.
Letting $\pi \colon \mathcal H \to  B(\mathcal H)$ be the multiplier representation that maps $f \in \mathcal H$ to the multiplication operator $(\pi f) \colon g \mapsto f g$, one readily verifies that the operator norm, $\lVert f\rVert_\text{op} := \lVert \pi f \rVert_{B(\mathcal H)}$, generates a coarser topology than the Hilbert space norm on $\mathcal H$, and satisfies
\begin{displaymath}
	\lVert fg\rVert_\text{op} \leq \lVert f\rVert_\text{op} \lVert g\rVert_\text{op}.
\end{displaymath}
Thus, $\mathcal H$ has the structure of a Banach algebra with respect to pointwise function multiplication.
When $k$ is a real-valued kernel we may further equip $\mathcal H$ with pointwise complex conjugation $^*\colon f \mapsto \bar f$ as an isometric involution, making it a Banach $^*$-algebra.
When the constant function $\bm 1\colon X \to \mathbb C$ with $\bm 1(x) = 1$ lies in $\mathcal H$, then $\mathcal H$ is a unital algebra with unit $\bm 1$.
Furthermore, the Hilbert space norm and operator norms become equivalent since $\lVert f \rVert_\mathcal H \leq \lVert \bm 1 \rVert_{\mathcal H} \cdot \lVert f \rVert_\text{op}$, i.e.,
\begin{displaymath}
	\frac{1}{\lVert \bm 1 \rVert_\mathcal H} \lVert f \rVert_\mathcal H \leq \lVert f \rVert_\text{op} \leq \lVert \Delta \rVert_{B(\mathcal H)} \lVert f \rVert_{\mathcal H}.
\end{displaymath}

In addition to having Banach algebra structure, an RKHA $\mathcal H$ is a cocommutative, coassociative coalgebra with $\Delta \colon \mathcal H \to \mathcal H \otimes \mathcal H$ as its comultiplication operator.
In particular, coassociativity of $\Delta$, i.e., $(\Delta \otimes \text{Id}) \circ \Delta = (\text{Id} \otimes \Delta) \circ \Delta$, is an important property that follows from associativity of multiplication (i.e., $\Delta^*$).
This property allows us to amplify $\Delta$ from $\mathcal H$ to the tensor product spaces $\mathcal H^{\otimes(n+1)}$, $n \in \mathbb N$, by defining $\Delta_n \colon \mathcal H \to \mathcal H^{\otimes (n+1)}$ as
\begin{equation*}
	\Delta_1 = \Delta, \quad \Delta_{n} = (\Delta \otimes \text{Id}^{\otimes (n-1)}) \Delta_{n-1} \quad \text{for } n>1.
\end{equation*}

Next, let $\sigma(\mathcal H)$ be the spectrum of an RKHA $\mathcal H$ as Banach algebra, i.e., the set of (automatically continuous) multiplicative linear functionals $\chi \colon \mathcal H \to \mathbb C$,
\begin{displaymath}
	\chi(fg) = (\chi f)(\chi g), \quad \forall f, g \in \mathcal H,
\end{displaymath}
equipped with the weak-$^*$ topology of $\mathcal H^* \supset \sigma(\mathcal H)$.
The pointwise evaluation functionals $\delta_x = \langle \varphi(x), \cdot \rangle_{\mathcal H}$ are elements of the spectrum for every $x \in X$, inducing a map $\hat\varphi\colon X \to \sigma(\mathcal H)$ by $\hat \varphi(x) = \delta_x$.
If the feature map $\varphi$ is injective, then so is $\hat\varphi$ and the spectrum $\sigma(\mathcal H)$ contains a copy of $X$ as a subset.

A useful property of some unital RKHAs is bijectivity of the map $\hat\varphi$.
This property implies (e.g., \cite{DasGiannakis23}*{Corollary~7}) that the spectrum $\sigma_{\mathcal H}(f)$ of each element $f \in \mathcal H$ (i.e., the set of complex numbers $z$ such that $f-z$ does not have a multiplicative inverse in $\mathcal H$) is equal to its range, $\sigma_{\mathcal H}(f) = f(X)$.
This enables in turn the use of holomorphic functional calculus of functions in $\mathcal H$,
\begin{displaymath}
	a(f) = \frac{1}{2\pi i} \int_\Gamma \frac{a(z)}{z-f} \, dz,
\end{displaymath}
where $a\colon D \to \mathbb C$ is any holomorphic function on a domain $D \subseteq \mathbb C$ that contains $\sigma_{\mathcal H}(f)$ and $\Gamma$ is an appropriate union of Jordan curves that encircle $\sigma_{\mathcal H}(f)$.
In particular, for every $\xi \in \mathcal H$ satisfying $\xi(x) \geq \varepsilon > 0$ for all $x \in X$, the $n$-th root $\xi^{1/n}$ lies in $\mathcal H$ for every $n \in \mathbb N$ by the holomorphic functional calculus.
The well-definition of such roots as elements of $\mathcal H$ plays a key role in the schemes discussed in \cref{sec:overview_fock_space}.

\subsection{Examples of RKHAs}

The main class of examples considered here are built on compact abelian groups, $G$.
Let $\hat G$ be the Pontryagin dual of $\hat G$, and $\mu$, $\hat{\mu}$ the Haar measures on $G$, $\hat G$, respectively, with $\mu$ normalized to a probability measure.
We use $\mathcal F \colon L^1(G) \to C_0(\hat{G})$ and $\hat{\mathcal F} \colon L^1(\hat{G}) \to C_0(G)$ to denote the Fourier transforms, and $*\colon L^1(\hat G) \times L^1(\hat G) \to L^1(\hat G)$ to denote convolution on the dual group:
$$(\mathcal F f)(\gamma) = \int_G f(x)\gamma(-x)\, d\mu(x), \quad (\hat{\mathcal F}\hat{f})(x) = \int_{\hat{G}} \hat{f}(\gamma)\gamma(x)\,d\hat{\mu}(x),$$
$$(\hat f * \hat g)(\gamma) = \int_{\hat G} \hat f(\gamma') g(\gamma-\gamma')\, d\hat\mu(\gamma'). $$
Note that, by compactness of $G$, the dual group $\hat G$ has a discrete topology and $\hat \mu$ is a counting measure.
See, e.g., \cite{Rudin17} for further details on Fourier analysis on groups.

By Bochner's theorem, for every positive function $\lambda \in L^1(\hat G)$ the translation-invariant function $k \colon G \times G \to \mathbb C$ defined as $k(x, x') = (\hat{\mathcal F} \lambda)(x - x')$ is a continuous, positive-definite kernel with an associated RKHS $\mathcal H$ of continuous functions.
Since $G$ is compact, the space $\mathcal H$ can be characterized in terms of a decay condition on Fourier coefficients,
\begin{displaymath}
	\mathcal H = \left\{ f \in C(G): \sum_{\gamma\in\hat G} \frac{\lvert \mathcal F f(\gamma)\rvert^2}{\lambda(\gamma)} < \infty \right\}.
\end{displaymath}
Moreover, by Mercer's theorem, the reproducing kernel $k$ can be expressed in terms of the uniformly convergent series 
\begin{displaymath}
	k(x, x') = \sum_{\gamma \in \hat G} \lambda(\gamma)\overline{\gamma(x)}\gamma(x'), \quad \forall x, x' \in G,
\end{displaymath}
and
\begin{equation}
	\label{eq:rkha_basis}
	\{ \psi_\gamma := \sqrt{\lambda(\gamma)}\gamma \}_{\lambda(\gamma)\neq 0}
\end{equation}
is an orthonormal basis of $\mathcal H$.
If $\lambda$ is a strictly positive function, then $\mathcal H$ is a dense subspace of $C(G)$ (equivalently, $k$ is a so-called universal kernel \cite{SriperumbudurEtAl11}) and the corresponding feature map $\varphi\colon G \to \mathcal H$ is injective.

\begin{defn}
	A strictly positive, summable function $\lambda \colon \hat G \to \mathbb R_{>0}$ is said to be \emph{subconvolutive} if there exists $C>0$ such that $\lambda * \lambda \leq C \lambda$ pointwise on $\hat G$.
\end{defn}

Subconvolutive functions have a rich history of study in the context of Beurling convolution algebras, e.g., \cite{Feichtinger79,Grochenig07,Kaniuth09}.
In \cite{GiannakisMontgomery25}, the RKHS $\mathcal H$ induced by a subconvolutive function $\lambda$ is shown to be an RKHA and $\Delta$ is easily diagonalizable by
$$\Delta\psi_\gamma = \sum_{\alpha+\beta=\gamma} \sqrt{\frac{\lambda(\alpha)\lambda(\beta)}{\lambda(\gamma)}} \psi_\alpha \otimes \psi_\beta, \quad \gamma \in \hat{G}.$$
Standard examples for $G = \mathbb T^N$, $\hat G \cong \mathbb Z^N$, include the families of functions $\{ \lambda_\tau \colon \mathbb Z^N \to \mathbb R_{>0} \}_{\tau>0}$ with subexponential decay, where
\begin{equation}
	\label{eq:lambda_subexp}
	\lambda_\tau(j) = \prod_{i=1}^N e^{-\tau \lvert j_i \rvert^p}
\end{equation}
for fixed $p \in (0,1)$.
For each such $p$, we have a one-parameter family of unital RKHAs $\mathcal H_\tau$ generated by a Markovian family of kernels $k_\tau$, i.e.,
\begin{equation*}
	k_\tau(x, x') \geq 0, \quad \int_G k_\tau(x, \cdot) \, d\mu = 1, \quad \forall x, x' \in X;
\end{equation*}
see \cite{DasGiannakis23}.
Moreover, in this setting the elements $\gamma$ of the dual group $\hat G$ are standard Fourier functions, $\gamma_j(x)=e^{i j \cdot x}$ for some $j \in \mathbb Z^N$ using the parameterization $x \in [0, 2\pi)^N$ for the $N$-torus.

\begin{rk}
	As an example involving non-compact abelian groups, $G = \mathbb R^N$, consider the following weights on $\hat G \cong \mathbb R^N$,
	\begin{equation}
		\label{eq:weights_rn}
		\lambda^{-1}(x) = e^{\tau \lvert x\rvert^p} (1+\lvert x\rvert)^s(\ln(e+\lvert x\rvert))^t.
	\end{equation}
	These weights are subconvolutive if either $0 <p<1$, $\tau >0$, $s,t \in \mathbb R$, or $p \in \{0,1\}$ and $s >N$ \cite{Feichtinger79,Grochenig07}.
\end{rk}

\begin{defn}
	A function $\lambda\colon \hat G \to \mathbb R_{>0}$ is said to satisfy the \emph{Gelfand--Raikov--Shilov (GRS)} condition if
	\begin{displaymath}
		\lim_{n \to \infty} \lambda(n\gamma)^{1/n} = 1, \quad \forall \gamma \in \hat{G},
	\end{displaymath}
	and the \emph{Beurling--Domar (BD)} condition if
	\begin{displaymath}
		\sum_{n=1}^\infty \dfrac{\ln(\lambda^{-1}(n\gamma))}{n^2} < \infty, \quad \forall \gamma \in \hat{G}.
	\end{displaymath}
\end{defn}

Functions in the families \eqref{eq:lambda_subexp} and~\eqref{eq:weights_rn} for $\tau>0$, $p \in (0,1)$, $s, t \geq 0$ satisfy both the GRS and BD conditions (in addition to being subconvolutive), which has two implications for $\mathcal H_\tau$.
First, the GRS condition implies that the associated map $\hat \varphi_\tau\colon G \to \sigma(\mathcal H_\tau)$ is a homeomorphism \cite{DasEtAl23,GiannakisMontgomery25}.
Second, as shown in \cite{GiannakisMontgomery25}, the BD condition implies that $\mathcal H_\tau$ contains functions with support contained in arbitrary compact sets with nonempty interior.
These conditions play a vital role for the construction of RKHAs on generic compact Hausdorff subspaces $X$ of $\mathbb T^N$ or $\mathbb R^N$ (e.g., attractors of dynamical systems).
We may restrict to $\mathcal H_\tau(X) := \overline{\text{span} \left\lbrace k_{\tau,x} \mid x \in X \right \rbrace}$ which is also an RKHA with comultiplication $\Delta|_{\mathcal H_\tau(X)}$.
More importantly, the spectrum of $\mathcal H_\tau(X)$ as a Banach algebra is $X$, hence, $\sigma_{\mathcal H_\tau(X)}(f) = f(X)$.
This enables the use of the holomorphic functional calculus of functions in $\mathcal H_\tau(X)$ as outlined in \cref{sec:rkha_defs}.
Note that existence of compactly supported functions on $X \subset \mathbb R^N$ automatically implies that $\mathcal H_\tau(X)$ contains non-analytic elements.
For examples of analytic extensions when the GRS condition is violated see \cite{DasEtAl23}*{section~3}. 

\section{Fock Space Amplification}
\label{sec:overview_fock_space}
Amplification is a common technique in operator algebras and representation theory.
Classical examples relevant to operator methods for dynamical systems are the Stinespring dilation theorem, which factorizes a completely positive map between $C^*$-algebras (including a quantum channel; see \cref{sec:qmda}) through a representation (``lifting'') of the input $C^*$-algebra on an auxiliary Hilbert space, and the Sz.-Nagy dilation theorem that lifts a contraction on a Hilbert space to a unitary on a larger Hilbert space \cite{Paulsen03}.
In this section, we employ Fock space constructions to embed spectrally regularized approximations, $U^t_\tau$, of the Koopman operator acting on RKHAs $\mathcal H_\tau$ in a higher dimensional space in a structure-preserving way.
We will see in \cref{Leibniz Rule} that this amplification will recover the Leibniz rule.
In quantum field theory, Fock spaces are used to model many-body quantum systems through a framework known as second quantization \cite{Lehmann04}.
It is therefore natural to view the methods described below as embeddings of classical dynamics into field-theoretic systems.

\subsection{Fock space constructions}
\label{sec:fock_space_construction}

We will consider two different amplifications from Fock space constructions.
Since the tensor algebra $T^\otimes(\mathcal H_\tau) := \mathbb C \oplus \mathcal H_\tau \oplus \mathcal H_\tau^{\otimes 2} \oplus \ldots$ is the free algebra on $\mathcal H_\tau$, we can represent multiplicative operators as a tensor product of operators and the algebra structure of $\mathcal H_\tau$ can be recovered by a quotient.
Let $F(\mathcal H_\tau) = \overline{T^\otimes(\mathcal H_\tau)}$ denote the full Fock space which is closed under the inner product
\begin{equation}
	\label{eq:fock_innerp}
	\begin{gathered}
		\langle a, b \rangle_{F(\mathcal H_\tau)} = \bar a b, \quad a, b \in \mathbb C, \\
		\langle a, f \rangle_{F(\mathcal H_\tau)} = 0, \quad a \in \mathbb C, \quad f \in \mathcal H_\tau, \\
		\langle f_1 \otimes \cdots \otimes f_n, g_1 \otimes \cdots \otimes g_n \rangle_{F(\mathcal H_\tau)} = \prod_{i=1}^n \langle f_i, g_i \rangle_{\mathcal H_\tau}, \quad f_i, g_i \in \mathcal H_\tau.
	\end{gathered}
\end{equation}
We may also define the weighted symmetric Fock space $\fk$ which is constructed as the closure of the symmetric tensor algebra $T^\vee(\mathcal H_\tau) := \mathbb C \oplus \mathcal H_\tau \oplus \mathcal H_\tau^{\vee 2} \oplus \ldots$ with respect to the inner product satisfying (cf.\ \eqref{eq:fock_innerp})
\begin{equation*}
	\begin{gathered}
		\langle a, b \rangle_\fk = \bar a b, \quad a, b \in \mathbb C, \\
		\langle a, f \rangle_\fk = 0, \quad a \in \mathbb C, \quad f \in \mathcal H_\tau, \\
		\langle f_1 \vee \cdots \vee f_n, g_1 \vee \cdots \vee g_n \rangle_{\fk} = \frac{w^2(n)}{n!^2} \sum_{\sigma,\sigma'\in S_n} \prod_{i=1}^n \langle f_{\sigma(i)}, g_{\sigma'(i)} \rangle_{\mathcal H_\tau}, \quad f_i, g_i \in \mathcal H_\tau,
	\end{gathered}
\end{equation*}
for a strictly positive weight function $w \colon \mathbb N_0 \to \mathbb R_{>0}$.
Here, $\vee$ denotes the symmetric tensor product, defined as the average
\begin{displaymath}
	f_1 \vee \dots \vee f_n = \frac{1}{n!} \sum_{\sigma \in S_n} f_{\sigma(1)} \otimes \dots \otimes f_{\sigma(n)}, \quad f_i \in \mathcal H_\tau,
\end{displaymath}
over the $n$-element permutation group $S_n$, and $\mathcal H_\tau^{\vee n}$ is the closed subspace of $\mathcal H_\tau^{\otimes n}$ consisting of symmetric tensors.
The map $f_1 \otimes \cdots \otimes f_n \mapsto f_1 \vee \cdots \vee f_n$ defines, by linear extension, the orthogonal projection from $\mathcal H_\tau^{\otimes n}$ to $\mathcal H_\tau^{\vee n}$.  We will use $\Omega \equiv 1 \in \mathbb C \subset \fk$ to denote the ``vacuum'' vector of the Fock space.
By convention, we will always choose $w$ such that $w(0) = 1$.

In \cite{GiannakisEtAl25}, it is shown that if $w^{-2}$ is summable and subconvolutive,
\begin{equation*}
	w^{-2} \in \ell^1(\mathbb N_0), \quad w^{-2} * w^{-2}(n) \leq C w^{-2}(n),
\end{equation*}
$\fk$ becomes a unital Banach algebra with respect to the symmetric tensor product for a norm $\vertiii{\cdot}_{\fk}$ equivalent to the Hilbert space norm,
\begin{displaymath}
	\vertiii{f \vee g}_{\fk} \leq \vertiii{f}_{\fk} \vertiii{g}_{\fk}, \quad \forall f, g \in \fk,
\end{displaymath}
and with $\Omega$ as the unit.
Moreover, associated with $\fk$ is a coproduct, i.e., a bounded operator $\Delta \colon \fk \to \fk \otimes \fk$ such that
\begin{displaymath}
	\Delta^*(f \otimes g) = f \vee g.
\end{displaymath}
Among many possible constructions, in this paper we use weights from the subexponential family
\begin{displaymath}
	w(n) = e^{\sigma n^p}, \quad \sigma > 0, \quad p \in (0, 1).
\end{displaymath}
Recall that this family of weights was also used to build the RKHAs $\mathcal H_\tau$ in~\cref{sec:rkha_defs}, though, in general, the weights used to build $\mathcal H_\tau$ and $F_w(\mathcal H_\tau)$ are independent.

The Banach algebra structure for the weighted Fock space is particularly useful since (as we will see in \cref{Leibniz Rule}) the spectrum of this Banach algebra gives a collection of spacial realizations of the approximate Koopman dynamics.
Towards this goal, we may characterize the spectrum of $\fk$ (i.e., the set of nonzero multiplicative functionals $\chi\colon \fk \to \mathbb C$) as the set
\begin{displaymath}
	\sigma(\fk) = \left\{ \chi = \langle \xi, \cdot \rangle_{\fk}: \xi = \sum_{n=0}^\infty w^{-2}(n) \eta^{\vee n}: \eta \in (\mathcal H_\tau)_{R_w} \right\} \subset \fk^*,
\end{displaymath}
where $R_w$ is the radius of convergence of the series $\sum_{n=1}^\infty w^{-2}(n) z^n$, $z \in \mathbb C$.
Since $w \in \ell^2(\mathbb N)$, we have $ R_w \geq 1 $ and the set of admissible vectors $\eta$ in the definition above includes the unit ball of $\mathcal H_\tau$.
The Riesz representation theorem provides a dual realization of the spectrum called the cospectrum,
\begin{displaymath}
	\sigma_\text{co}(\fk) = \left\{ \xi \in \fk : \langle \xi, \cdot \rangle_{\fk} \in \sigma(\fk) \right\} \subset \fk.
\end{displaymath}
Equivalently, we have that $\sigma_\text{co}(\fk)$ is the subset of $\fk$ consisting of elements $\xi$ such that $\Delta \xi = \xi \otimes \xi$.
We equip $\sigma(\fk)$ and $\sigma_\text{co}(\fk)$ with the weak-$^*$ topology on $\fk^*$ and the weak topology on $\fk$, respectively.
With these topologies, they become compact Hausdorff spaces.

The weighted Fock space $\fk$ may be realized as an RKHA of continuous functions on $C(\sigma(\fk))$ via the Gelfand map $\Gamma(\xi)(\chi) = \chi(\xi)$; see~\cite{GiannakisMontgomery25}.
In addition, $\fk$ has an associated feature map $\tilde\varphi_\tau\colon X \to \fk$, where
\begin{equation*}
	\tilde\varphi_\tau (x) = \sum_{n=0}^\infty \frac{w^{-2}(n)}{\varpi_\tau^n} \varphi_\tau(x)^{\vee n},
\end{equation*}
$\varphi_\tau\colon X \to \mathcal H_\tau$ is the canonical feature map associated with $\mathcal H_\tau$ (see~\cref{sec:rkha_defs}), and $\varpi_\tau$ is a constant chosen such that $\varpi_\tau \geq \sup_{x \in X} \lVert \varphi_\tau(x)\rVert_{\mathcal H_\tau}$.
The range of $\tilde\varphi_\tau$ then lies in the cospectrum $\sigma_\text{co}(\fk)$, which implies that $\hat \varphi_\tau \colon X \to \sigma(\fk)$ with
\begin{equation*}
	\hat\varphi_\tau(x) = \langle \tilde\varphi_\tau(x), \cdot\rangle_{\fk}
\end{equation*}
is a well-defined map of state space into the spectrum of the weighted Fock space.
We will denote the image of state space under this map by $\hat X_\tau = \hat\varphi_\tau(X)$.

\subsection{Spectral approximation}\label{sec:spectral_approximation}

We would like to leverage the functional analytic properties of the RKHAs $\mathcal H_\tau$ and associated Fock spaces to build embeddings of unitary Koopman dynamics on $H = L^2(\mu)$.
However, for a general dynamical flow $\Phi^t \colon X \to X$ and RKHS $\mathcal H_\tau$ on $X$ the Koopman operator $f \mapsto f \circ \Phi^t$ need not be well-defined as a linear (let alone unitary) map from $\mathcal H_\tau$ to itself.
Even if $U^t\colon \mathcal H_\tau \to \mathcal H_\tau$ were a well-defined unitary, intuition from the $L^2$ setting suggests that for sufficiently complex dynamics it would be a non-diagonalizable operator with non-trivial continuous spectrum, hampering computations that rely on its eigendecomposition (see, e.g., \cref{sec:qc}).
To overcome this limitation, we regularize the skew-adjoint generator $V\colon D(V) \to H$ to build an approximating family of skew-adjoint operators $W_\tau \colon D(W_\tau) \to \mathcal H_\tau$ on the RKHAs $\mathcal H_\tau$.
These operators are diagonalizable and generate unitaries $U^t_\tau = e^{t W_\tau}$ that approximate the Koopman operator $U^t \colon H \to H$ in a suitable sense.

With the notation of \cref{sec:rkha_defs}, let $K_\tau\colon H \to \mathcal H_\tau$ be the (Hilbert--Schmidt) integral operator associated with the reproducing kernel of $H_\tau$, $K_\tau f = \int_X k_\tau(\cdot, x) f(x) \, d\mu(x)$.
By standard RKHS results, $K_\tau^* \colon \mathcal H_\tau \to H$ implements the restriction map from pointwise-defined functions in $\mathcal H_\tau$ to their corresponding equivalence classes in $H$.
Defining $G_\tau \colon H \to H$ as $G_\tau = K_\tau^* K_\tau$, we will assume that the RKHSs $\mathcal H_\tau$ are built such that $\{ G_\tau \}_{\tau>0} \cup \Id$ forms a semigroup of Markov operators with $\ran G_\tau \subset D(V)$.
In particular, we have that the family $\{ G_\tau\}_{\tau>0}$ is an approximate identity on $H$, $\lim_{\tau\to 0^+} \lVert G_\tau f - f \rVert_{H} = 0$ for every $f \in H$.
It can also be shown that $K_\tau$ admits the polar decomposition $K_\tau = T_\tau  G_\tau$, where $T_\tau \colon H \to \mathcal H_\tau$ is a partial isometry that satisfies $T_\tau f = K_{\tau/2} f$.
See, e.g., \cite{GiannakisEtAl24,GiannakisEtAl25} for further details on building kernel families with these properties.

Next, we will make use of the following notions of convergence of skew-adjoint operators; e.g., \cite{Chatelin11,Oliveira09}.

\begin{defn} Let $A\colon D(A) \to \mathbb H$ be a skew-adjoint operator on a Hilbert space $\mathbb H$ and $A_\tau \colon D(A_\tau) \to \mathbb H$ a family of skew-adjoint operators indexed by $\tau$, with resolvents $R_z(A) = (zI - A)^{-1}$ and $R_z(A_\tau) = (zI - A_\tau)^{-1}$ respectively, for a complex number $z$ in the resolvent sets of $A$ and $A_\tau$.
	\begin{enumerate}
		\item The family $A_\tau$ is said to converge in \emph{strong resolvent sense} to $A$ as $\tau\to 0^+$ if for some (and thus, every) $z \in \mathbb C \setminus i \mathbb R$ the resolvents $R_z(A_\tau)$ converge strongly to $R_z(A)$; that is, $\lim_{\tau\to 0^+}R_z(A_\tau) f = R_z(A) f$ for every $f\in \mathbb H$.
		\item The family $A_\tau$ is said to converge in \emph{strong dynamical sense} to $A$ as $\tau\to 0^+$ if for every $t \in \mathbb R$, the unitary operators $e^{tA_\tau}$ converge strongly to $e^{tA}$; that is, $\lim_{\tau\to0^+} e^{tA_\tau} f = e^{tA} f$ for every $f\in \mathbb H$.
	\end{enumerate}
\end{defn}

It can be shown, e.g., \cite{Oliveira09}*{Proposition~10.1.8}, that strong resolvent convergence and strong dynamical convergence are equivalent notions.
For our purposes, this implies that if a family of skew-adjoint operators $V_\tau$ on $H$ converges to the Koopman generator $V$ in strong resolvent sense, the unitary evolution groups $e^{t V_\tau}$ generated by these operators consistently approximate the Koopman group $U^t = e^{t V}$ generated by $V$.
Strong resolvent convergence and strong dynamical convergence also imply a form of strong convergence of spectral measures, e.g., \cite{DasEtAl21}*{Proposition~13}.

With these definitions, our Fock space amplification schemes are based on a family of approximating operators $V_\tau\colon D(V_\tau) \to H$ with the following properties:

\begin{enumerate}[label=(V\arabic*)]
	\item \label[prty]{prty:V1} $V_\tau$ is skew-adjoint.
	\item \label[prty]{prty:V2} $V_\tau$ is real, i.e., $\overline{V_\tau f} = V_\tau \bar f$ for all $f \in D(V_\tau)$.
	\item \label[prty]{prty:V3} $V_\tau$ is diagonalizable.
	\item \label[prty]{prty:V4} $V_\tau$ has a simple eigenvalue at 0 with $\bm 1$ as a corresponding eigenfunction.
	\item \label[prty]{prty:V5} $D(V_\tau)$ includes $\ran G_{\tau/2}$ as a subspace.
	\item \label[prty]{prty:V6} As $\tau \to 0^+$, $V_\tau$ converges in strong resolvent sense, and thus in strong dynamical sense, to $V$.
\end{enumerate}

Methods for building families of operators $V_\tau$ satisfying \crefrange{prty:V1}{prty:V6} can be found in \cite{DasEtAl21,GiannakisValva24,GiannakisValva25}, but any other technique with similar properties can be employed in what follows.

Let $X_\mu = \supp(\mu)$ be the (compact) support of the invariant measure $\mu$ and define the skew-adjoint operators $W_\tau \colon D(W_\tau) \to \mathcal H_\tau$ on the dense domains $D(W_\tau) = T_\tau(D(V_\tau)) \oplus \mathcal H_\tau(X_\mu)^\perp $  as $W_\tau = T_\tau V_\tau T_\tau^*$.
Note that the well-definition of $W_\tau$ depends on \cref{prty:V5}.
By \cref{prty:V3}, for every $\tau>0$, $W_\tau$ is a diagonalizable operator.
Moreover, by \cref{prty:V4}, the restriction of $W_\tau$ on $\mathcal H_\tau(X_\mu)$ admits the eigendecomposition
\begin{displaymath}
	W_\tau \zeta_{j,\tau} = i \omega_{j,\tau} \zeta_{j,\tau}, \quad j \in \mathbb N_0, \quad \omega_{0,\tau} = 0, \quad \omega_{2j,\tau} = -\omega_{2j-1,\tau}, \quad \zeta_{2j,\tau} = \overline{\zeta_{2j-1,\tau}},
\end{displaymath}
with $\omega_{j,\tau} \in \mathbb R$ and $\{ \zeta_{j,\tau}\}_{j \in \mathbb N_0}$ forming an orthonormal basis of $\mathcal H_\tau$.
We interpret $\zeta_{j,\tau}$ in the above as approximate Koopman eigenfunctions and $\omega_{j,\tau}$ as their corresponding eigenfrequencies.
Note that, unlike Koopman eigenfunctions in $L^2$, the $\zeta_{j,\tau}$ are pointwise-defined functions on $X$ with regularity controlled by the reproducing kernel $k_\tau$.

By skew-adjointness, $V_\tau$ and $W_\tau$ generate strongly continuous, one-parameter unitary groups on $H$ and $\mathcal H_\tau$, respectively.
These unitaries, $e^{t V_\tau}$ and $e^{t W_\tau}$, can be used interchangeably to approximate the Koopman evolution of observables in $H$, as follows.

\begin{lem}[\cite{GiannakisEtAl24}*{Lemma~8}]
	\label{lem:koopman_approx}
	For every $f \in H$, $t \in \mathbb R$, and $\tau>0$, we have $K_\tau^* e^{t W_\tau} K_\tau f= G_{\tau/2}e^{tV_\tau} G_{\tau/2}f$.
	Moreover, the following hold.
	\begin{enumerate}
		\item For every $f \in H$, $\lim_{\tau\to 0^+} K_\tau^* e^{t W_\tau} K_\tau f = U^t f$.
		\item For every $f \in \mathcal H_1$, $\lim_{\tau\to 0^+} K_\tau^* e^{t W_\tau} f = U^t \iota f$.
	\end{enumerate}
\end{lem}

Despite being dynamically consistent as $\tau \to 0^+$, the approximate generators $V_\tau$ and $W_\tau$ have an important structural difference compared to the Koopman generator $V$, namely, they generally fail to satisfy a Leibniz rule (cf.~\eqref{eq:leibniz}).
In particular, despite being well-defined elements of the RKHA $\mathcal H_\tau$, pointwise products $\zeta_{i,\tau} \zeta_{j,\tau}$ of eigenfunctions of $W_\tau$ are not necessarily eigenfunctions.
Yet, the fact that $W_\tau$ is ``close'' (in the sense of strong resolvent convergence) to the generator $V$ (which does obey Leibniz), suggests that products such as $\zeta_{i,\tau} \zeta_{j,\tau}$ should contain useful dynamical information.
In the following subsections, we will pass to the free-algebra setting of the Fock spaces described in \cref{sec:fock_space_construction} in order to capture this information.

\subsection{Recovering the Leibniz rule}\label{Leibniz Rule}

In an effort to utilize more information from a limited set of approximate eigenfunctions $\zeta_{j,\tau}$, and block-diagonalize an amplified approximate Koopman operator, we will use the Leibniz rule and nontrivial liftings in the Fock space (see \cites{GiannakisEtAl24,GiannakisEtAl25}).

\subsubsection{Amplification to the full Fock space}

The scheme of \cite{GiannakisEtAl24} lifts the regularized generator $W_\tau$ to a skew-adjoint operator $\tilde W_\tau \colon D(\tilde W_\tau) \to F(\mathcal H_\tau) $, defined by linear extension of
\begin{equation}
	\label{eq:lift_w}
	\begin{aligned}
		\tilde W_\tau(f_1 \otimes f_2 \otimes \cdots \otimes f_n)
		 & = (W_\tau f_1) \otimes f_2 \otimes \cdots \otimes f_n \\
		 & + f_1 \otimes (W_\tau f_2) \otimes \cdots \otimes f_n \\
		 & + \ldots                                              \\
		 & + f_1 \otimes f_2 \otimes \cdots \otimes (W_\tau f_n)
	\end{aligned}
\end{equation}
for $f_1, \ldots, f_n \in D(W_\tau)$.
By construction, $\tilde W_\tau$ satisfies the Leibniz rule
\begin{displaymath}
	\tilde W_\tau(f \otimes g) = (\tilde W_\tau f) \otimes g + f \otimes (\tilde W_\tau g)
\end{displaymath}
with respect to the tensor product for all $f, g \in D(\tilde W_\tau)$ such that the left- and right-hand sides of the above equation are well-defined.
As a result, we have:
\begin{enumerate}
	\item The point spectrum $\sigma_p(\tilde W_\tau)$ is a union of abelian subgroups of $i \mathbb R$ generated by $\sigma_p(W_\tau)$.
	\item $\tilde W_\tau$ generates a 1-parameter group of unitary operators $\tilde U^t_\tau = e^{t \tilde W_\tau}$, $ t \in \mathbb R$, that act multiplicatively with respect to the tensor product,
	      \begin{displaymath}
		      \tilde U^t_\tau(f \otimes g) = (\tilde U^t_\tau f) \otimes (\tilde U^t_\tau g) , \quad \forall, f, g \in F(\mathcal H_\tau).
	      \end{displaymath}
\end{enumerate}

Putting together the above, it follows that for every vector $q \in L^2(\mu)$ with a representative $\xi \in \mathcal H_{\tau_0}$ for some $\tau_0>0$ and a multiplicative decomposition of the form $\xi = \xi_1 \ldots \xi_n$ for some $\xi_1, \ldots, \xi_n \in \mathcal H_{\tau_0}$, we have
\begin{equation*}
	\Delta_{n-1}^* \tilde U^t_\tau \xi = \Delta_{n-1}^* \left( \bigotimes_{i=1}^n U^t_\tau \xi_i \right) = \prod_{i=1}^n U^t_\tau \xi_i \xrightarrow{\tau \to 0^+} U^t q
\end{equation*}
in $L^2(\mu)$ norm, where $\Delta_{n-1}^*\colon \mathcal H_\tau^{\otimes n} \to \mathcal H_\tau$ is the $n$-fold multiplication
\begin{equation*}
	\Delta_{n-1}^* = \Delta^* \circ (\Id_{\mathcal H_\tau} \otimes \Delta_{n-2}^*), \quad \Delta^*_1 = \Delta^*.
\end{equation*}
Intuitively, the Fock space $F(\mathcal H_\tau)$ generated by the RKHA $\mathcal H_\tau$ allows one to ``distribute'' the Koopman evolution of observables and states over tensor products in the Fock space for potentially arbitrarily high grading $n \in \mathbb N$.

\subsubsection{Amplification to the weighted symmetric Fock space}

The same analysis above applies to the weighted Fock space approach.
However, the Banach algebra structure of $F_w(\mathcal H_\tau)$ provides a new dynamical interpretation.
Fix eigenfunctions of $W_\tau$, $\{\zeta_{i,\tau}\}_{i=0}^\infty$ as in \cref{sec:spectral_approximation}.
We will see that the amplified approximate Koopman operator is given by precomposition by an irrational rotation on any one of a parameterized family of possibly infinite dimensional tori.

Let $\mathbb A_w$ be the subset of $\ell^2(\mathbb N_0)$ consisting of vectors $a = (a_j)_j$ with norm $\lVert a\rVert_{\ell^2(\mathbb N)} \leq R_w$ and non-negative elements $a_j$.
For each $a \in \mathbb A_w$ and each sequence $z = (z_j)_j \in \ell^\infty(\mathbb N)$ with unimodular elements $z_j$, define the vectors
\begin{equation*}
	\xi_{\tau,a,z} = \sum_{n=0}^\infty w^{-2}(n) \left(a_0 +\sum_{j=1}^\infty a_j z_j \zeta_{j,\tau} \right)^{\vee n} \in \sigma_\text{co}(\fk),
\end{equation*}
and the subsets $\mathbb T_{\tau, a}$ of the spectrum $\sigma(\fk)$ as
\begin{displaymath}
	\mathbb T_{\tau, a} = \left\{ \chi_{\tau,a,z} \equiv \langle \xi_{\tau,a,z}, \cdot\rangle_{\fk}: z = (z_j)_j \in \ell^\infty(\mathbb N), \; z_j \in \mathbb T^1 \subset \mathbb C \right\}.
\end{displaymath}
Each set $\mathbb T_{\tau, a}$ has the topology of a torus of dimension equal to the number of nonzero elements of $(a_1, a_2, \ldots)$.
In what follows, $\mathbb S_\tau = \bigcup_{a \in \mathbb A_w} \mathbb T_{\tau, a} \subset \sigma(\fk)$ will be the (disjoint) union of these tori.

Let us define a lifting $\tilde W_\tau\colon D(\tilde W_\tau) \to \fk$ similarly to~\eqref{eq:lift_w}, but with $\otimes$ replaced by the symmetric tensor product $\vee$.
Then the approximate Koopman operators $\tilde U^t_\tau = e^{t W_\tau}$ induce a flow $R_\tau^t \colon \sigma(\fk) \to \sigma(\fk)$, $t \in \mathbb R$, where $R^t_\tau(\chi) = \chi \circ \tilde U^t_\tau$.
Furthermore, each torus $\mathbb T_{\tau, a}$ is an invariant set under $R_\tau^t$.
On these sets, $R^t_\tau$ takes the form of a rotation system generated  by the eigenfrequencies $\omega_{j,\tau}$:
\begin{displaymath}
	R_\tau^t(\chi_{\tau,a,z}) = \chi_{\tau,a, z^{(t)}},
\end{displaymath}
where $z = (z_j)_{j \in \mathbb N}$ and $z^{(t)} = (e^{-i \omega_{j,t}} z_j)_{j \in \mathbb N}$.
Equivalently, we have
\begin{displaymath}
	R_\tau^t(\langle \xi_{\tau,a,z}, \cdot\rangle_\fk) = \langle \tilde U^{-t}_\tau\xi_{\tau,a,z}, \cdot\rangle_\fk,
\end{displaymath}
so the vector $\xi_{\tau,a,z}$ evolves under the adjoint (``Perron--Frobenius'') operators $\tilde U^{t*}_\tau = \tilde U^{-t}_\tau$.

The rotation system $R_\tau^t$ constitutes a topological model of the regularized Koopman dynamics $U^t_\tau$ as a rotation system on the spectrum of the weighted Fock space $\fk$.
This is non-trivial since $U^t_\tau$ is not a composition operator induced by a flow on the original state space $X$.
Since $\hat X_\tau$ is a subset of $\mathbb S_\tau \subset \sigma(\fk)$ (and $\mathbb S_\tau$ is invariant under $R^t_\tau$), the union of tori $\mathbb S_\tau$ provides a common topological setting for studying the dynamical system associated with the regularized Koopman operators $U^t_\tau$ (represented by $R_\tau^t$) in relation to the original dynamical system $\Phi^t$ (represented by its image $\hat X_\tau$ under the map $\hat\varphi_\tau$).

One important question is whether each of these tori and rotation systems remember the approximate Koopman dynamics.
If you choose the weights $a \in \mathbb A_w$ to be strictly positive, then the values of a function, $\xi \in F_w(\mathcal H_\tau)$, on the torus $\mathbb T_{\tau,a}$ uniquely determines $\xi$ (see \cite{GiannakisEtAl25}).
This means the rotation system on $\mathbb T_{\tau,a}$ uniquely determines the amplified approximate Koopman operator, $\tilde{U}_\tau^t$.
We will see in \cref{sec:qc} how to implement rotation systems on tori as  quantum circuits.

\subsection{Embeddings of observables and state vectors}

To utilize the Leibniz rule in the Fock space, we must lift both the observable and the state vector to the Fock space.
The two schemes we will discuss here come from two different algebra homomorphisms,
\begin{gather*}
	\Delta_{n-1}^* \colon \mathcal H_\tau^{\otimes n} \to \mathcal H_\tau, \\
	|_X \colon F_w(\mathcal H_\tau) \to F_w(\mathcal H_\tau)(X),
\end{gather*}
where $|_X$ denotes the restriction map.
In both approaches we make use of a higher dimensional embedding of observables or states that map down to the desired observable or state via $\Delta_n^*$ or $|_X$.
We will refer to these higher dimensional embeddings as liftings (e.g. $\tilde{f} \in F_w(\mathcal H_\tau)$ is a lifting of $f \in F_w(\mathcal H_\tau)(X)$ if $\tilde{f}|_X = f$).

\subsubsection{Tensor network approach}

We will begin with the scheme \cite{GiannakisEtAl24} that can be implemented as a tree tensor network.
Given a state vector $\xi \in \mathcal H_\tau$ we will choose $\tilde{\xi} \in \mathcal H_\tau^{\otimes n}$ such that $\Delta_{n-1}^*\tilde{\xi} = \xi$.
Moreover, given an observable $f \in L^\infty(\mu)$ and regularization parameters $\sigma,\tau>0$, we choose the mapping
\begin{equation}
	\label{eq:tensornet_obs}
	f \mapsto A_{f,\sigma,\tau,n} = K_\sigma^{\otimes n} \Delta_{n-1} K_\tau M_f K_\tau^* \Delta_{n-1}^* (K_\sigma^*)^{\otimes n} \in B(\mathcal H_\tau^{\otimes n})
\end{equation}
into quantum observables acting on the tensor product space $\mathcal H_\tau^{\otimes n}$.
Note that~\eqref{eq:tensornet_obs} is positivity preserving similarly to $\pi\colon L^\infty(\mu) \to B(H)$ employed in \cref{sec:qmda}, but it is not an algebra homomorphism on $L^\infty(\mu)$.
Disregarding the smoothing by the kernel operator $K_\sigma^{\otimes n}$, we see how finite multiplicativity is enforced in the quantum expectation of $A_{f,\sigma,\tau,n}$ with respect to the pure quantum state $\rho \propto \langle (\xi^{1/n})^{\otimes n}, \cdot \rangle_{\mathcal H_\tau^{\otimes n}}(\xi^{1/n})^{\otimes n}$,
\begin{equation*}
	\left\langle \Delta_{n-1} K_\tau M_f K_\tau^* \Delta_{n-1}^* (U_\tau^t \xi^{1/n})^{\otimes n}, (U_\tau^t \xi^{1/n})^{\otimes n}\right\rangle_{\mathcal H_\tau^{\otimes n}} =
	\left\langle K_\tau M_f K_\tau^* (U_\tau^t \xi^{1/n})^{n}, (U_\tau^t \xi^{1/n})^{n}\right\rangle_{\mathcal H_\tau^{\otimes n}}.
\end{equation*}
In particular, at time $t=0$ we also recover the initial expectation.
See \cite{GiannakisEtAl24}*{Figure~4}  for a schematic illustrating how the computations above can be organized in a tensor network.
In \cref{prop:pointwise} below, we will renormalize by the expectation of the constant observable since the operators $\Delta_{n-1}$ are not unit norm, and can lead to a rescaled expectation.

The states used for the tensor network scheme may be built so as to approximate statistical expectations of observables in $L^\infty(\mu)$ with respect to probability densities in $L^1(\mu)$ (as in~\cref{sec:qmda}), as well as pointwise evaluation of continuous observables $f \in C(X)$.
As an illustration of the latter, let the state space $X$ be the $N$-dimensional torus, and let $p_{\bm\mu,\bm\kappa} \in C^\infty(\mathbb T^N) \cap \mathcal H_\tau$ denote the von Mises distribution with locality parameters $\bm \mu = (\mu_1, \ldots, \mu_N) \in [0, 2\pi)^N$ and concentration parameters $\bm \kappa = (\kappa_1, \ldots, \kappa_N) \in \mathbb R_+^N$,
			$$p_{\bm \mu, \bm \kappa}= \prod_{i=1}^N p_{\mu_i, \kappa_i}\qquad p_{\mu,\kappa}(\theta) = \frac{e^{\kappa\cos(\theta-\mu)}}{I_0(\kappa)}.$$
			For a probability distribution $(w_i)_{i=1}^\infty$ set
			\begin{gather*}
				\xi = F_{\kappa,\tau}(x) = p_{x,\bm \kappa} / \lVert p_{x,\bm \kappa}\rVert_{\mathcal H_\tau}, \quad \bm \kappa = (\kappa, \ldots, \kappa) \in \mathbb R^N_+, \\
				\eta_\tau = w_1 \frac{\xi}{\lVert \xi\rVert_{\mathcal H_\tau}} + w_2 \frac{\xi^{1/2} \otimes \xi^{1/2}}{\lVert \xi^{1/2}\rVert_{\mathcal H_\tau}^2} + w_3 \frac{\xi^{1/3} \otimes \xi^{1/3} \otimes \xi^{1/3}}{\lVert \xi^{1/3}\rVert_{\mathcal H_\tau}^3} + \ldots,
			\end{gather*}
			and
			$$\tilde \Xi_{\kappa,\tau}(x) = \langle \eta_\tau, \cdot \rangle_{F(\mathcal H_\tau)} \eta_\tau.$$
			In the above, we interpret $\tilde\Xi_{\kappa,\tau}\colon X \to B_1(F(H_\tau))$ as a ``quantum feature map" (whereas $F_{\kappa,\tau}\colon X \to \mathcal H_\tau$ is an RKHA-valued classical feature map).
			We let $\tilde{\mathcal P}^t_\tau \colon B_1(F_w(\mathcal H_\tau)) \to B_1(F_w(\mathcal H_\tau))$ be the quantum transfer operator induced from $\tilde U^t_\tau$, defined analogously to $\mathcal P\colon B_1(H) \to B_1(H)$ from \cref{sec:qm_overview}.
			We then have the following pointwise approximation result for the Koopman evolution of continuous observables.

			\begin{prop}
				\label{prop:pointwise}
				Fix $f \in C(\mathbb T^N)$, and with the assumptions of \cite{GiannakisEtAl24}, define $\tilde f^{(t)}_{\kappa,\sigma,\tau,n} \in C(\mathbb T^N)$ where $\tilde f^{(t)}_{\kappa,\sigma,\tau,n}(x)$ is obtained from a normalized expectation of $A_{f,\sigma,\tau,n}$ with respect to the quantum state $\tilde \Xi_{\kappa,\tau}(x)$, i.e.,
				\begin{displaymath}
					\tilde f^{(t)}_{\kappa,\sigma,\tau,n}(x) = \frac{\mathbb E_{\tilde{\mathcal P}^t_\tau(\tilde\Xi_{\kappa,\tau}(x))} A_{f, \sigma, \tau,n}}{\mathbb E_{\tilde{\mathcal P}^t_\tau(\tilde\Xi_{\kappa,\tau}(x))} A_{\bm 1, \sigma, \tau,n}}.
				\end{displaymath}
				Then, for every $x \in X$, $\lim_{\kappa\to\infty}\lim_{\sigma\to 0^+}\lim_{\tau\to 0^+} \tilde f^{(t)}_{\kappa,\sigma,\tau,n}(x) = U^t f(x)$.
			\end{prop}

			\subsubsection{Second quantization approach}
			\label{sec:second_quantization}

			For the second quantization scheme, given $f \in F_w(\mathcal H_\tau)(X) =: \tilde{\mathcal H}_\tau$ we will pick an observable $\tilde{f} \in F_w(\mathcal H_\tau)$ such that $\tilde{f}|_X =f$.
			The state vectors on the other hand are provided by a feature map from $X$ into $\sigma(F_w(\mathcal H_\tau))$.
			This allows a significant amount of freedom to choose such a lifting.
			The ones chosen below were picked mainly out of convenience for numerical experiments.

			We will use the identification of $\fk$ with the RKHA $\hfk \subset C(\sigma(\fk))$ of functions from \cref{sec:fock_space_construction}.
			Integral operators may then be used to represent $f \in H$ by elements of the RKHA $\hfk$, whose restrictions on finite-dimensional tori $\mathbb T_{\tau, a} \subset \mathbb S_\tau$ are polynomials of arbitrarily large degree $m \in \mathbb N$ of the coordinates $z_j$.

			Let $\kappa: X \times X \to \mathbb R_{>0}$ be a strictly positive, bounded, continuous kernel function such that
			\begin{equation}
				\label{eq:kappa_decay}
				\kappa(x, y) < \kappa(x, x), \quad \forall x, y \in X: \, x \neq y.
			\end{equation}
			For example, given a metric $d \colon X \times X \to \mathbb R$ that metrizes the topology of $X$, a prototypical kernel satisfying \eqref{eq:kappa_decay} is the radial Gaussian kernel,
			\begin{equation*}
				\kappa(x,y) = \exp\left( - \frac{d^2(x,y)}{\varepsilon^2}\right), \quad \varepsilon > 0.
			\end{equation*}
			For any such kernel $\kappa$ and $\tau > 0$, define the smoothed kernel $\kappa_\tau \colon X \times X \to \mathbb R_{> 0}$, where $\kappa_\tau(\cdot, y) = K_\tau \iota \kappa(\cdot, y) \in \mathcal H_\tau$.
			Moreover, for $m \in \mathbb N$ define the integral operators $\mathcal K_{m,\tau} \colon H \to \fk$ and $\hat{\mathcal K}_{m,\tau} \colon H \to \hfk$, where
			\begin{displaymath}
				\mathcal K_{m,\tau} f = \int_X \kappa_\tau^{\vee m}(\cdot, y) f(y) \, d\mu(y), \quad \hat{\mathcal K}_{m,\tau} f = \Gamma \mathcal K_{m,\tau} f = \int_X (\Gamma \kappa_\tau(\cdot, y))^m f(y) \, d\mu(y).
			\end{displaymath}
			Well-definition of these operators is verified in \cite{GiannakisEtAl25}.
			Note that for every $x \in X$ the pointwise power $y \mapsto (\kappa_\tau(x, y))^m$ lies in the RKHA $\tilde{\mathcal H}_\tau \subset C(X)$ with reproducing kernel
			\begin{displaymath}
				\tilde k_\tau(x, y) = \sum_{n=0}^\infty \frac{w^{-2}(n)}{\varpi_\tau^n} k_\tau(x, y)^n,
			\end{displaymath}
			where $\varpi_\tau$ is a normalization function chosen so as to ensure that $\tilde k_\tau(x,\cdot)$ lies in the cospectrum of $F_w(\mathcal H_\tau)$ for every $x \in X$.
			Moreover, $\kappa_\tau^{\vee m}(\cdot, y) \in \fk$ lies above $\kappa_\tau(\cdot, y)^m \in \tilde{\mathcal H}_\tau$,
			\begin{displaymath}
				\tilde \pi(\kappa_\tau^{\vee m}(\cdot, y)) = \kappa_\tau(\cdot, y)^m,
			\end{displaymath}
			where $\tilde{\pi} \colon F_w(\mathcal H_\tau) \to F_w(\mathcal H_\tau)(X)$ is the unique orthogonal projection onto $\fk(X)$.

			We map observables $f \in H$ to elements of the RKHA $\hfk$ by means of the integral operators $\hat{\mathcal K}_{m,\tau}$; specifically, $\hat g_{m,\tau} := \hat{\mathcal K}_{m,\tau} \iota f$.
			We will also employ $\hat h_{m,\tau} := \hat{\mathcal K}_{m,\tau} 1_X$ for normalization purposes. 
			The functions $\hat g_{m,\tau}$ and $\hat h_{m,\tau}$ evolve unitarily under the action of the Koopman operator $\hat U^t_\tau$ to
			\begin{displaymath}
				\hat g^{(t)}_{m,\tau} := \hat U^t_\tau \hat g_{m,\tau}, \quad \hat h^{(t)}_{m,\tau} := \hat U^t_\tau \hat h_{m,\tau},
			\end{displaymath}
			respectively.

			In order to render an approximation of the true Koopman evolution $U^t f$ using $\hat g^{(t)}_{m,\tau}$ and $\hat h^{(t)}_{m,\tau}$, we pull back these functions to state space $X$ by means of a feature map.
			Setting $\sigma > 0$ and $ \tau \in (0, \sigma/2]$, we let $\varphi_\sigma \colon X \to \mathcal H_\sigma$ be the canonical feature map into $\mathcal H_\sigma$,
\begin{displaymath}
	\varphi_\sigma(x) = k_\sigma(x, \cdot),
\end{displaymath}
and $\varphi^{(\mu)}_\sigma$ its projection onto the subspace $\mathcal H_\sigma(X_\mu) \subseteq \mathcal H_\sigma$,
\begin{displaymath}
	\varphi^{(\mu)}_\sigma = \proj_{\mathcal H_\sigma(X_\mu)} \circ \varphi_\sigma.
\end{displaymath}
We then define $\tilde\varphi^{(\mu)}_\sigma \colon X \to F_w(\mathcal H_\sigma)$
\begin{displaymath}
	\tilde\varphi^{(\mu)}_\sigma (x) = \sum_{n=0}^\infty \frac{w^{-2}(n)}{\varpi_\sigma^{2n}} \varphi^{(\mu)}_\sigma(x)^{\vee n}, \quad \varpi_\sigma = \sup_{x \in X} \lVert \varphi_\sigma(x)\rVert_{\mathcal H_\sigma}.
\end{displaymath}
Since $\sigma > \tau$, we have $ F_w(\mathcal H_\sigma) \subset \fk $, so we can view $\tilde\varphi^{(\mu)}_\sigma$ as a feature map into $\fk$ that is based on the canonical feature vectors associated with $\mathcal H_\sigma(X_\mu) \subseteq \mathcal H_\sigma \subseteq \mathcal H_\tau$.
In particular, for every $\sigma \geq \tau$ we have $\tilde\varphi^{(\mu)}_\sigma(x) \in \sigma_\text{co}(\fk)$ since $\varpi_\sigma = \lVert \varphi_\sigma(x)\rVert_{\mathcal H_\sigma} \geq \lVert \varphi_\sigma(x)\rVert_{\mathcal H_\tau}$ and $\lVert \varphi^{(\mu)}_\sigma(x)\rVert_{\mathcal H_\sigma} \leq \lVert \varphi_\sigma(x)\rVert_{\mathcal H_\sigma}$.

We therefore obtain a feature map $\hat\varphi^{(\mu)}_{\sigma,\tau} \colon X \to \sigma(\fk)$ mapping into the spectrum of $\fk$ via 
\begin{displaymath}
	\hat\varphi^{(\mu)}_{\sigma,\tau}(x) = \langle \tilde\varphi^{(\mu)}_\sigma(x), \cdot\rangle_{\fk}.
\end{displaymath}
By construction, the image $\hat X^{(\mu)}_{\sigma,\tau} = \hat\varphi_{\sigma,\tau}(X)$ of state space under this feature map lies in in the union of tori $\mathbb S_\tau \subset \sigma(\fk)$, so we can dynamically evolve each point $\hat\varphi_{\sigma,\tau}(x)$ using the rotation system $R^t_\tau$ restricted to $\mathbb S_\tau$.

\begin{rk}
	A reason for building the two-parameter family of feature maps $\hat\varphi^{(\mu)}_{\sigma,\tau}$ based on $\varphi^{(\mu)}_\sigma$ (as opposed to, say, the canonical feature maps $\varphi_\tau$ of $\mathcal H_\tau$) is to control the regularity of feature vectors as elements of $H$ when taking $\tau \to 0^+$ limits associated with Koopman operator approximation.
	Indeed, it follows from the semigroup structure of $\{ G_\tau \}_{\tau >0}$ (see \cref{sec:spectral_approximation}) that for every $\tau \in (0, \sigma/2]$ and $x \in X$,
	\begin{equation*}
		\varphi^{(\mu)}_\sigma(x) = K_\tau G_{\frac{\sigma}{2}-\tau} q_\sigma(x),
	\end{equation*}
	for a continuous function $q_\sigma \colon X \to H$.
	Then strong convergence (\cref{lem:koopman_approx}) implies
	\begin{equation*}
		\lim_{\tau\to 0^+} \left\lVert K_\tau^* U^t_\tau \varphi^{(\mu)}_\sigma(x) - U^t \iota \varphi^{(\mu)}_\sigma(x)\right\rVert_H = 0, \quad \forall x \in X.
	\end{equation*}
	One further finds that
	\begin{equation*}
		q_\sigma(x) = K^*_{\sigma/2} \varphi_{\sigma/2}(x),
	\end{equation*}
	which allows to relate $H$ inner products with $q_\sigma(x)$ to pointwise evaluation in the RKHS $\mathcal H_{\sigma/2}$,
	\begin{displaymath}
		\langle q_\sigma(x), f \rangle_H = \langle \varphi_{\sigma/2}, K_{\sigma/2}f\rangle_{\mathcal H_{\sigma/2}} = (K_{\sigma/2} f)(x), \quad \forall f \in H, \quad \forall x \in X.
	\end{displaymath}
\end{rk}

Next, it can be shown that sufficiently small $\tau$, $\hat h^{(t)}\rvert_{\hat X^{(\mu)}_{\sigma,\tau}}$ is bounded away from zero and thus has a multiplicative inverse in $C(\hat X^{(\mu)}_{\sigma,\tau})$.
As a result,
\begin{equation*}
	\hat f_{m,\tau}^{(t)} := \frac{\hat g^{(t)}_{m,\tau}\rvert_{\hat X^{(\mu)}_{\sigma,\tau}}}{\hat h^{(t)}_{m,\tau}\rvert_{\hat X^{(\mu)}_{\sigma,\tau}}}, \quad f^{(t)}_{m,\sigma,\tau} := \hat f_{m,\tau}^{(t)} \circ \hat \varphi^{(\mu)}_{\sigma,\tau} \rvert_X \equiv \frac{\hat g^{(t)}_{m,\tau} \circ \hat\varphi^{(\mu)}_{\sigma,\tau}\rvert_X}{\hat h^{(t)}_{m,\tau} \circ \hat\varphi^{(\mu)}_{\sigma,\tau}\rvert_X}
\end{equation*}
are well-defined continuous functions on $\hat X^{(\mu)}_{\sigma,\tau}$ and $X$, respectively.
As a result, we have:
\begin{prop}[\cite{GiannakisEtAl25}*{Theorem~6}]
	With notation and assumptions as above, $f^{(t)}_{m,\sigma,\tau} \to U^t f$ $\mu$-a.e.\ according to the limit
	\begin{displaymath}
		\lim_{m\to\infty} \lim_{\sigma \to 0^+} \lim_{\tau \to 0^+}\left\lVert \iota(f^{(t)}_{m, \sigma, \tau} - U^t f)\right\rVert_H = 0.
	\end{displaymath}
\end{prop}

\section{Quantum simulation of systems with pure point spectra}%
\label{sec:qc}

In this section we restrict attention to  measure-preserving, ergodic flows with pure point spectra; i.e., continuous-time measure-preserving ergodic systems for which the union of eigenspaces of the Koopman operator is dense in $L^p(\mu)$, $p \in [1, \infty)$.
These systems are examples of highly structured dynamics that in many ways is antithetical to weak-mixing.

As mentioned in~\cref{sec:intro}, by the Halmos--von Neumann theorem \cite{HalmosVonNeumann42}, every ergodic pure-point-spectrum system is isomorphic in a measure-theoretic sense to a rotation system on a compact abelian group.
In continuous time, the point spectrum $\sigma_p(V)$ of the generator is an additive subgroup of the imaginary line, giving, by the spectral mapping theorem, the point spectra of the Koopman operators as a group homomorphism, $\sigma_p(V) \ni i \omega \mapsto e^{i\omega t} \in \sigma_p(U^t)$.
As we will see below, for this restricted class of systems the group structure of the spectrum of the generator enables a tensor product factorization of the associated unitary Koopman group that admits a highly efficient gate-based implementation on a quantum computer.

As a concrete example, we consider an ergodic rotation $\Phi^t\colon X \to X$ on the torus $X= \mathbb T^d$,
\begin{displaymath}
	\Phi^t(\theta) = \theta + \alpha \cdot \theta \mod 2 \pi, \quad t \in \mathbb R,
\end{displaymath}
where $\alpha = (\alpha_1, \ldots, \alpha_d) \in \mathbb R^d$ are rationally-independent frequency parameters.
This system is a canonical representative in the measure-theoretic isomorphism class of pure-point-spectrum, continuous-time ergodic systems with spectra generated by $d$ basic frequencies.
Specifically, we have
\begin{equation}
	\label{eq:rot_spec}
	\sigma_p(V) = \{ i(j_1 \alpha_1 + \cdots + j_d \alpha_d): j_1, \ldots, j_d \in \mathbb Z \},
\end{equation}
so the point spectrum of the generator is isomorphic to the Pontryagin dual $\widehat{\mathbb T^d} \simeq \mathbb Z^d$ of the state space of the dynamics by rational independence of $\alpha_1, \ldots, \alpha_d$.

\subsection{RKHA embedding of dynamics}

Let $\mathcal H$ be an RKHA of continuously differentiable functions on $\mathbb T^d$ built using a subconvolutive weight function on $\mathbb Z^d$; e.g., the subexponential weights $\lambda^{-1/2}$ from~\eqref{eq:lambda_subexp}.
The Koopman operators $U^t\colon \mathcal H \to \mathcal H$, $U^t f = f \circ \Phi^t$ are well-defined unitaries on $\mathcal H$ and form a strongly continuous group generated by a skew-adjoint operator $V\colon D(V) \to \mathcal H$, $D(V) \subset \mathcal H$, with point spectrum given by~\eqref{eq:rot_spec} as in the $L^p$ case.
The corresponding orthonormal eigenvectors can be chosen as the scaled characters $\psi_j$ from~\eqref{eq:rkha_basis}; that is, the generator $V$ on the RKHA $\mathcal H$ admits the diagonalization
\begin{equation*}
	V \psi_j = i \omega_j \psi_j, \quad \omega_j = j_1 \alpha_1 + \cdots + j_d \alpha_d.
\end{equation*}
Here, the group structure of $\sigma_p(V)$ is a manifestation of the fact that $V$ obeys an analog of the Leibniz rule~\eqref{eq:leibniz} on the RKHA $\mathcal H$.
In the $L^2$ setting, \cite{TerElstLemanczyk17} have shown that satisfying~\eqref{eq:leibniz} on the algebra $L^\infty(\mu) \cap D(V)$ is a necessary and sufficient condition for a skew-adjoint operator $V\colon D(V) \to L^2(\mu)$ to be the generator of a unitary Koopman group.
Given that the Leibniz rule is a fundamental structural property of vector fields generating continuous-time deterministic dynamical systems, we conjecture that analogous result holds for unitary groups of composition operators on RKHAs.

In the context of quantum computing, RKHAs usefully combine the ability to represent classical observables $f \in \mathcal H$ by bounded multiplication operators $\pi f \in B(\mathcal H)$ with the ability to implement pointwise evaluation at $x \in X$ as a quantum mechanical expectation.
Specifically, for $x\in X$ we define the vector state $\mathbb E_x\in S_*(B(\mathcal H))$ as
\begin{equation*}
	\mathbb E_x A := \frac{\langle k_x, A k_x \rangle_{\mathcal H}}{k(x, x)}.
\end{equation*}
Then, for every observable $f \in \mathcal H$ and evolution time $t \in \mathbb R$ we have
\begin{equation*}
	f(\Phi^t(x)) = \mathbb E_x(\mathcal U^t M_f),
\end{equation*}
where $\mathcal U^t : B(\mathcal H) \to B(\mathcal H)$ is the conjugation operator $\mathcal U^t A = U^t A U^{t*}$ defined analogously to the $L^2$ case (see \cref{sec:intro}).
Equivalently, we can consider the evolution of normal states under the dual operator $\mathcal P^t\colon S_*(B(\mathcal H)) \to S_*(B(\mathcal H))$, $\mathcal P^t \rho = U^{t*} \rho U^t$.
Defining the quantum feature map $\Xi\colon X \to S_*(B(\mathcal H))$ that sends $x \in X$ to the rank-1 density operator associated with $\mathbb E_x$, $\rho_x = \langle k_x, \cdot \rangle_{\mathcal H} k_x / k(x, x)$, we have that $\Xi$ intertwines \emph{pointwise} classical dynamics on $X$ with quantum dynamics on $S_*(B(\mathcal H))$, i.e.,
\begin{equation*}
	\mathcal P^t \circ \Xi = \Xi \circ \Phi^t.
\end{equation*}

We can lift this correspondence to the level of Borel probability measures by recalling that the set of Borel probability measures on a compact Hausdorff space $X$ may be identified with the state space of the $C^*$-algebra of continuous functions on $X$, $C(X)$, equipped with the supremum norm.
The state space $X$ embeds weak-$^*$ continuously into $S(C(X))$ under the map $\delta\colon X \to S(C(X))$ that sends points to their corresponding Dirac measures, $\delta(x) = \delta_x$.
This map intertwines state space dynamics under $\Phi^t$ with the dynamics of probability measures under the pushforward map $\Phi^t_* \colon S(C(X)) \to S(C(X))$, $\Phi^t_* \circ \delta = \delta \circ \Phi^t$.

Next, using the quantum feature map $\Xi$, we define a map $Q\colon S(C(X)) \to S_*(B(\mathcal H))$ from Borel probability measures on $X$ to normal quantum states on $\mathcal H$ via the operator-valued (Bochner) integral $Q(\nu) = \int_X \Xi(x) \, d\nu(x)$.
The map $Q$ is injective (by universality of the reproducing kernel of $\mathcal H$) and weak-$^*$ to weak continuous and it intertwines classical statistical evolution under $\Phi^t_*$ with quantum evolution under $\mathcal P^t$,
\begin{equation*}
	\mathcal P^t \circ Q = Q \circ \Phi^t_*.
\end{equation*}
Moreover, we have $\Xi = Q \circ \delta$ by construction.

The relationships outlined above can be summarized through the following commuting diagram:
\begin{equation}
	\label{eq:commut_rkha}
	\begin{tikzcd}
		X \ar[r,"\Phi^t"] \ar[d,"\delta",swap] & X \ar[d,"\delta"] \\
		S(C(X)) \ar[r,"\Phi^t_*"] \ar[d,"Q",swap] & S(C(X)) \ar[d,"Q"] \\
		S_*(B(\mathcal H)) \ar[r,"\mathcal P^t"] & S_*(B(\mathcal H))
	\end{tikzcd}.
\end{equation}
This diagram is more general than the commuting diagram~\eqref{eq:commut_transfer} for measure-preserving dynamics in $L^p$, in the sense that~\eqref{eq:commut_rkha} captures the evolution of arbitrary Borel probability measures (as opposed to $\mu$ absolutely continuous measures in~\eqref{eq:commut_transfer}), including the evolution of individual states in $X$ through their representation as Dirac measures.
Of course, it should be kept in mind that in order to achieve this level of generality in the dynamical relationships between classical and quantum levels of description in this subsection we have significantly restricted the class of systems under study to pure point spectrum systems.

Before closing this subsection we note that (unlike $L^2$ spaces) multiplication operators $M_f$ by real-valued elements $f \in \mathcal H$ are generally not self-adjoint, and thus do not constitute quantum observables that can be physically measured (e.g., at the output stage of a quantum circuit).
To overcome this issue we can represent $f \in \mathcal H$ by the self-adjoint operator $S_f = (M_f + M_f^*)/2 \in B(\mathcal H)$; one then verifies that $\mathbb E_x S_f = \mathbb E_x M_f$ whenever $f$ and the reproducing kernel of $\mathcal H$ are real-valued.

\subsection{Quantum simulation algorithm}

Let $\mathbb B$ be the abstract two-dimensional Hilbert space over the complex numbers with orthonormal basis $ \{ \ket 0, \ket 1 \}$.
We identify $\mathbb B$ with the quantum mechanical Hilbert space associated with a qubit.
A quantum computer with $n$ qubits is associated with the tensor product Hilbert space $\mathbb B_n = \mathbb B^{\otimes n}$ of dimension $2^n$.
We will use the notation $\ket b = \ket{b_0} \otimes \cdots \otimes \ket{b_n}$ with $b = (b_1, \ldots, b_n) \in \{ 0, 1 \}^n$ to denote the elements of the tensor product basis of $\mathbb B_n$ induced from the $ \{ \ket 0, \ket 1 \}$ basis of $\mathbb B$, known in this context as the quantum computational basis.
For later convenience, we introduce the Pauli $Z$ operator on $\mathbb B$, which is defined via the eigenvalue/eigenvector relations $Z \ket 0 = \ket 0$ and $Z \ket 1 = - \ket 1$.
For a detailed exposition of quantum computing we refer the reader to \cite{NielsenChuang10}.

We build a quantum algorithm for approximating the quantum dynamics on $\mathcal H$ described above and, by virtue of \eqref{eq:commut_rkha}, the classical dynamics on $\mathbb T^d$ by defining a family of projections $\hat \Pi_q\colon \mathcal H \to \mathbb B_n$ with $q \in \mathbb N$ and $n = d(q+1)$ that carry along the generator $V$ on $\mathcal H$ to skew-adjoint operators $\hat V_q \colon  \mathbb B_n \to \mathbb B_n$, $\hat V_q = \hat \Pi_q V \hat \Pi_q$.
These projections are chosen such that the corresponding unitary evolution operators $\hat U^t_q = e^{t \hat V_q}$ on $\mathbb B_n$ can be implemented by a quantum algorithm whose running time is bounded by a constant that is independent of $n$, using quantum circuits with $O(n)$ quantum logic gates.
Moreover, as $q$ increases, the quantum expectations from the quantum computational system recover the evolution of observables from the original system.
Specifically, we have
\begin{equation*}
	\lim_{q\to\infty} \tr((\hat{\mathcal P}^t_q \hat\Xi_q(x)) (\hat\pi_q f)) = \mathbb E_x(\mathcal U^t M_f)= f(\Phi^t(x))
\end{equation*}
for all $f \in \mathcal H$, $x \in X$, and $t \in \mathbb R$, where $\hat\Xi_q\colon X \to S_*(B(\mathbb B_n))$ and $\hat\pi_q \colon \mathcal H \to B(\mathbb B_n)$ are projected quantum feature maps and representations of observables, respectively,
\begin{displaymath}
	\hat\Xi_q(x) = \frac{\hat \Pi_q \Xi(x) \hat \Pi_q}{\tr(\hat \Pi_q \Xi(x) \hat \Pi_q)}, \quad \hat\pi_q f = \hat \Pi_q (\pi f) \hat \Pi_q,
\end{displaymath}
and $\hat{\mathcal P}^t_q\colon S_*(B(\mathbb B_n)) \to S_*(B(\mathbb B_n))$ is the quantum state evolution $\hat{\mathcal P}^t_q\rho = \hat U^{t*}_q \rho \hat U^t_q$ induced by $\hat U^t_q$.

To build $\hat \Pi_q$ we start from the $2^n$-dimensional subspaces $\mathcal H_q \subset \mathcal H$, defined as
\begin{displaymath}
	\mathcal H_q = \spn \{ \psi_j : \text{$j = (j_1, \ldots, j_d) \in J_q^d$} \}
\end{displaymath}
for $J_q = \{ -2^q,  \ldots, -1, 1, \ldots, 2^q \}$.
Note that 0 is excluded from the index set $J_q$.
To every $j_i \in J_q$ we assign a binary string $\beta_q(j_i) \in \{ 0, 1 \}^{q+1}$ given by the binary representation of the integer $\tilde j_i \in \{ 0, \ldots, 2^{q+1} -1 \}$ with
\begin{displaymath}
	\tilde j_i =
	\begin{cases}
		j_i + 2^q,    & j_i < 0, \\
		j_i + 2^q - 1 & j_i > 0,
	\end{cases}
\end{displaymath}
and to every $j \in J_q^d$ we assign $\bm \beta_q(j) \in \{ 0, 1 \}^{d(q+1)}$ by concatenation, $\bm \beta_q(j) = (\beta_q(j_1), \ldots, \beta_q(j_d))$.
We then define $\hat \Pi_q$ by linear extension of
\begin{displaymath}
	\hat \Pi_q \psi_j =
	\begin{cases}
		\ket{\bm \beta_q(j)}, & j \in J_q^d,      \\
		0,                    & \text{otherwise}.
	\end{cases}
\end{displaymath}
It can be shown \cite{GiannakisEtAl22} that the projected generator takes the form
\begin{equation}
	\label{eq:vdecomp}
	\begin{aligned}
		\hat V_q & = \hat v_{1} Z \otimes I \otimes I \otimes I \otimes \cdots \otimes I \\
		         & + \hat v_{2} I \otimes Z \otimes I \otimes I \otimes \cdots \otimes I \\
		         & + \hat v_{4} I \otimes I \otimes Z \otimes I \otimes \cdots \otimes I \\
		         & + \ldots                                                              \\
		         & + \hat v_{2^{n-1}} I \otimes \cdots \otimes I \otimes Z,
	\end{aligned}
\end{equation}
where $\hat v_1, \ldots, \hat v_{2^{n-1}}$ are complex valued coefficients determined by computing a Fourier--Walsh transform \cite{WelchEtAl14} applied to the spectrum $\sigma_p(\hat V_q)$ of the projected generator.

It follows from~\eqref{eq:vdecomp} that the unitary $\hat U^t_q = e^{t\hat V_q}$ admits the factorization
\begin{equation}
	\label{eq:hat_u_t_q}
	\hat U^t_q = \bigotimes_{i=0}^{n-1} e^{t \hat v_{2^i} Z},
\end{equation}
which can be implemented by phase rotations acting on $n$ qubits in parallel  with no cross-cubit communication.
The resulting quantum algorithm thus has constant running time and can be implemented using circuits of depth $O(n)$, as mentioned above.
It is important to note that the Leibniz rule~\eqref{eq:leibniz} is crucial to the efficiency of this algorithm.
In particular, the ability to decompose the generator in the form of~\eqref{eq:vdecomp}, rather than a more general Fourier--Walsh expansion, relies on the group structure of the point spectrum $\sigma_p(V)$ implied by the Leibniz rule.
It should also be noted that in order to perform a full computation the circuit implementing~\eqref{eq:hat_u_t_q} must be pre- and post-composed with appropriate state preparation and measurement stages.
The depth of the entire circuit is then dominated by a quantum Fourier transform at the measurement stage which is $O(n^2)$; see \cite{GiannakisEtAl22}.

\bigskip

\subsection*{Acknowledgments} The authors acknowledge support from the U.S.\ Department of Defense, Basic Research Office under Vannevar Bush Faculty Fellowship grant N00014-21-1-2946, the U.S.\ Office of Naval Research under MURI grant N00014-19-1-242, and the U.S.\ Department of Energy under grant DE-SC0025101.


\end{document}